\documentclass[11pt]{article}
\usepackage{amssymb}
\usepackage{amsmath}
\usepackage{graphicx}
\usepackage{acronym}
\usepackage{amsfonts}
 \usepackage{cite}
 \usepackage{epstopdf}

\newcommand{\qed}{\bull \medskip}

\newtheorem{theorem}{Theorem}
\newtheorem{corollary}{Corollary}
\newtheorem{proposition}{Proposition}
\newtheorem{lemma}{Lemma}
\newtheorem{definition}{Definition}
\newtheorem{algorithm}{Algorithm}
\newtheorem{example}{Example}
\def\be{\begin{example}}
\def\ee{\end{example}}
\def\bt{\begin{theorem}}
\def\et{\end{theorem}\bigskip}
\def\bl{\begin{Lemma}}
\def\el{\end{Lemma}\bigskip}
\def\ep{\end{Proposition}\bigskip}
\def\bp{\begin{Proposition}}
\def\bd{\begin{definition}}
\def\ed{\end{definition}}

\newcommand{\alglist}{
\begin{list}{Step 1}
{\setlength{\leftmargin}{1.1 in}\setlength{\labelwidth}{1.0 in}} }

\def\qed{\hfill {$ \Box $} \medskip}
\def\proof{\noindent\bf Proof. \hspace{1mm}\rm}

\newcommand{\x}{{\bf x}}
\newcommand{\y}{{\bf y}}

\newcommand{\T}{{\cal T}}

\newcommand{\A}{{\cal A}}
\newcommand{\B}{{\cal B}}

\begin{document}
\title{\bf  A Necessary and Sufficient Condition for Existence of a Positive Perron Vector}
 \author{Shenglong Hu\footnote{Department of Mathematics, School of Science, Tianjin University, Tianjin, 300072 China.
 Computational and Applied Mathematics Initiative, Department of Statistics, University of Chicago, Chicago, IL 60637-1514.
 His work is partially supported by National Science Foundation and National Science Foundation of China (Grant No. 11401428). Email: timhu@tju.edu.cn; tim.hu@galton.uchicago.edu.}\quad and Liqun Qi\footnote{Department of Applied
 Mathematics, The Hong Kong Polytechnic University, Hung Hom, Kowloon, Hong Kong.
His work was supported by the Hong Kong Research Grant Council
(Grant No. PolyU 502111, 501212, 501913 and 15302114). E-mail:
maqilq@polyu.edu.hk.}}

\date{\today}
\maketitle

{\large

\begin{abstract}
In 1907, Oskar Perron showed that a positive square matrix has a unique largest positive eigenvalue with a positive eigenvector.   This result was extended to irreducible nonnegative matrices by Geog Frobenius in 1912, and to irreducible nonnegative tensors and weakly irreducible nonnegative tensors recently.  This result is a fundamental result in matrix theory and has found wide applications in probability theory, internet search engines, spectral graph and hypergraph theory, etc.  In this paper, we give a necessary and sufficient condition for the existence of such a positive eigenvector, i.e., a positive Perron vector, for a nonnegative tensor.  We show that every nonnegative tensor has a canonical nonnegative partition form, from which we introduce strongly nonnegative tensors. A tensor is called strongly nonnegative, if the spectral radius of each genuine weakly irreducible block is equal to the spectral radius of the tensor, which is strictly larger than the spectral radius of any other block. We prove that a nonnegative tensor has a positive Perron vector if and only if it is strongly nonnegative.   The proof is nontrivial.  Numerical results for finding a positive Perron vector are reported.
\medskip


\end{abstract}

\medskip

\section{Introduction}

More than one century ago, in 1907, Oskar Perron showed that a positive square matrix has a unique largest positive eigenvalue with a positive eigenvector (Perron vector) \cite{P}.   This result was extended to irreducible nonnegative matrices by Geog Frobenius \cite{F} in 1912, and to irreducible nonnegative tensors by Chang, Pearson and Zhang \cite{cpz08} in 2008, and weakly irreducible nonnegative tensors by Friedland, Gaubert and Han \cite{FGH} in 2013.  This result is a fundamental result in matrix theory \cite{BP, hj85} and has found wide applications in probability theory \cite{BP, hj85}, the Google Page Rank \cite{Page}, spectral graph and hypergraph theory \cite{cd12, q14}, etc.

Perhaps the most important part of the Perron-Frobenius theorem, as well as a key feature ensuring Google Page Rank's success, is the assertion that
\begin{multline*}
\begin{matrix}\textit{The nonnegative Perron vector is unique for an irreducible nonnegative matrix and}\\
\textit{it is the positive Perron vector.}\end{matrix}
\end{multline*}
Thus, irreducibility is a sufficient condition to guarantee the existence of a positive Perron vector.
In this paper, we will study
\[
\textit{The necessary and sufficient condition for the existence of a positive Perron vector.}
\]
We will study this problem in a general setting for all nonnegative tensors of orders higher than or equal to two, which includes the nonnegative matrix case since matrices are second order tensors.    In the matrix case, this result is known \cite[Page 10-11]{Ro}.    However, in the case of order higher than two, the proof is nontrivial.   

The following theorem will be proved.
\begin{theorem}\nonumber
A nonnegative tensor has a positive Perron vector if and only if it is strongly nonnegative.
\end{theorem}
Strongly nonnegative tensors will be defined in Definition~\ref{def:strong}, whereas the matrix counterpart will be presented in Section~\ref{sec:matrix} for illustration.

\subsection{Strongly Nonnegative Matrix}\label{sec:matrix}
Given an $n\times n$ nonnegative matrix $A$, we can always partition $A$ (up to some permutations) into the following upper triangular block form
\begin{equation}\label{matrix}
A=\begin{bmatrix}A_1&A_{12}&\dots&\dots&\dots&\dots&\dots&A_{1r}\\ &A_2&A_{23}&\dots&\dots&\dots&\dots&A_{2r}\\ & &\ddots&\ddots &&&&\vdots\\ &&& A_s&A_{ss+1}&\dots&\dots&A_{sr}\\ &&&&A_{s+1}&0&\dots&0\\ &&&&&\ddots&&\vdots\\ &&&&&&\ddots\\&&&&&&& A_r
\end{bmatrix}
\end{equation}
such that:
\begin{enumerate}
\item Each diagonal block matrix $A_i$ is irreducible for $i\in[r]:=\{1,\dots,r\}$. Here we regard a scalar (zero or not) as a one dimensional irreducible matrix for convenience.
\item For each $i\in [s]$, at least one of the matrices $A_{ij}$ is not zero for $j=i+1,\dots,r$.
\end{enumerate}
Then the matrix $A$ is \textit{strongly nonnegative} if
\[
\rho(A_i)<\rho(A)\ \text{for all }i\in [s]\ \text{and }\rho(A_i)=\rho(A)\ \text{for all }i=s+1,\dots,r.
\]

\subsection{Outline}
In Section~\ref{sec:preliminary}, we will present some basic definitions and results on nonnegative tensors. In Section~\ref{sec:partition} we will first review the nonnegative tensor partition result from \cite{HHQ}, and then refine the partition by introducing \textit{genuine weakly irreducible principal sub-tensors}. In Section~\ref{sec:nnt}, we will give the necessary and sufficient condition for a nonnegative tensor possessing a positive spectral radius. This class of nonnegative tensors is called \textit{nontrivially nonnegative}. In Section~\ref{sec:positive-vector}, we will first introduce \textit{strongly nonnegative tensors}, and then prove that a tensor being strongly nonnegative is both necessary and sufficient for it to have a positive Perron vector. The proofs in both Sections~\ref{sec:nnt} and \ref{sec:positive-vector} are based on the nonnegative tensor partition developed in Section~\ref{sec:partition}. In Section~\ref{sec:algoas}, we will propose an algorithm to determine whether a nonnegative tensor being strongly nonnegative or not, and find a positive Perron vector when it is.

\section{Preliminaries}\label{sec:preliminary}
Eigenvalues of tensors, independently proposed by  Qi \cite{Qi} and Lim \cite{l05} in 2005, have been becoming active research topics in numerical multilinear algebra and beyond. We refer to \cite{cpz08,CQZ,cd12,DW2,DW,GLY,Ka,oo12,q14,YY2,YY} and references therein for some recent developments and applications. Among others, nonnegative tensors have wide investigations, see \cite{cpz08,FGH,HHQ,YY2,YY} for the Perron-Frobenius type theorems, and \cite{cd12,DQW,DW,q14,GLY} for some applications; and the survey \cite{CQZ} for more connections to other problems in hypergraph, quantum entanglement, higher order Markov chain, etc.

An $m$th order $n$-dimensional tensor $\A$ is a multi-way array $\A=(a_{i_1\dots i_m})$ indexed by $m$ indices $i_j$ for $j\in[m]$ with each $i_j$ being within the set $[n]$. Usually, the entries $a_{i_1\dots i_m}$ can be elements in any prefixed set $S$, not necessarily scalars; and the set of tensors of order $m$ and dimension $n$ with entries in $S$ is denoted by $T_{m,n}(S)$.
The space of tensors of order $m$ and dimension $n$ with entries in the field $\mathbb C$ of complex numbers is denoted simply by $T_{m,n}$.
In this article, we will focus on the case when $S$ is the nonnegative orthant $\mathbb R_+$. The interior of $\mathbb R_+$ is denoted by $\mathbb R_{++}$. A tensor $\A=(a_{i_1\dots i_m})$ with $a_{i_1\dots i_m}\in\mathbb R_+$ for all $i_j\in [n]$ and $j\in [m]$ is called a \textit{nonnegative tensor}.
Let $N_{m,n}\subset T_{m,n}$ be the set of all nonnegative tensors with order $m$ and dimension $n$. Therefore, $N_{2,n}$ represents the set of all nonnegative $n\times n$ matrices.

For any $j\in [n]$, let
\begin{equation}\label{eqn:index-j}
I(j):=\{(i_2,\dots,i_m)\in [n]^{m-1}: j\in\{i_2,\dots,i_m\}\}.
\end{equation}
Obviously, $I(j)$ depends on $m$.
For example, $I(j)=\{j\}$ when $m=2$. We omit this dependence notationally for simplicity, since $m$ is always clear from the content.
For each $\A\in N_{m,n}$,
we associate it a nonnegative matrix $M_{\A}=(m_{ij})\in N_{2,n}$ with
\[
m_{ij}=\sum_{(i_2,\dots,i_m)\in I(j)}a_{ii_2\dots i_m}.
\]
The matrix $M_{\A}$ is called the \textit{majorization matrix} of $\A$ (cf.\ \cite{HHQ,P10}).
\begin{definition}[Weakly Irreducible Nonnegative Tensor \cite{FGH,HHQ}]\label{def:weak}
A nonnegative tensor $\A\in N_{m,n}$ is called \textit{weakly irreducible} if the majorization $M_{\A}$ is irreducible. $\A$ is weakly reducible if it is not weakly irreducible.
\end{definition}
For convenience, tensors in $N_{m,1}=\mathbb R_+$ are always regarded as weakly irreducible. Note that the weak irreducibility for tensors in $N_{2,n}$ (i.e., matrices) reduces to the classical irreducibility for nonnegative matrices (cf.\ \cite{BP,hj85}).

For any nonnegative matrix $A\in N_{2,n}$, we associate it a directed graph $G=(V,E)$ as $V=\{1,\dots,n\}$ and
\[
(i,j)\in E\ \text{if and only if }a_{ij}>0.
\]
It is known that the irreducibility of the matrix $A$ is equivalent to the strong connectedness of the corresponding directed graph defined as above \cite{hj85}.

\begin{definition}[Eigenvalues and Eigenvectors \cite{Qi,l05}]\label{def:eigenvalue-eigenvector}
Let tensor $\mathcal A=(a_{i_1\ldots i_m})\in T_{m,n}$. A number $\lambda\in\mathbb C$ is called an \textit{eigenvalue} of $\mathcal A$, if there exists a vector $\mathbf x\in\mathbb C^n\setminus\{\mathbf 0\}$ which is called an \textit{eigenvector} such that
\begin{equation}\label{eigenvalue-equation}
\mathcal A\mathbf x^{m-1}=\lambda\mathbf x^{[m-1]},
\end{equation}
where $\mathbf x^{[m-1]}\in\mathbb C^n$ is an $n$-dimensional vector with its $i$-th component being $x_i^{m-1}$, and $\mathcal A\mathbf x^{m-1}\in\mathbb C^n$ with
\[
\big(\mathcal A\mathbf x^{m-1}\big)_i:=\sum_{i_2,\dots,i_m=1}^na_{ii_2\dots i_m}x_{i_2}\dots x_{i_m}\ \text{for all }i=1,\dots,n.
\]
\end{definition}
The number of eigenvalues of a tensor is always finite (cf. \cite{Qi}). The spectral radius $\rho(\A)$ of a tensor $\A$ is defined as
\[
\rho(\A):=\max\{|\lambda| : \lambda\ \text{is an eigenvalue of }\A\}.
\]
In general $\rho(\A)$ is not an eigenvalue of $\A$, however it is when $\A$ is nonnegative \cite{YY}. There are extensive research on Perron-Frobenius type theorems for nonnegative tensors, see \cite{cpz08,FGH,YY,HHQ,YY2,CQZ} and references therein.

In the following, we summarize the Perron-Frobenius theorem for nonnegative tensors which will be used in this article.
\begin{theorem}[Perron-Frobenius theorem \cite{YY,FGH,hj85}]\label{thm:pf}
Suppose that $\A\in N_{m,n}\setminus\{0\}$.  Then, the following results hold.
\begin{itemize}
\item [(i)] $\rho(\A)$ is an eigenvalue of $\A$ with an eigenvector in $\mathbb R^n_+$.
\item[(ii)] If $\A$ has a positive eigenvector $\y\in\mathbb R^n_{++}$ associated with an eigenvalue $\lambda$, then $\lambda=\rho(\A)$ and
$$
\min_{\mathbf x\in\mathbb R^n_{++}}\max_{1\leq i\leq
n}\frac{\left(\A\x^{m-1}\right)_i}{x_i^{m-1}}=\rho(\A)=\max_{\x\in\mathbb R^n_{++}}\min_{1\leq
i\leq n}\frac{\left(\A\x^{m-1}\right)_i}{x_i^{m-1}}.
$$
\item[(iii)] If $\A$ is weakly irreducible, then $\A$ has a positive eigenpair $(\lambda=\rho(\A),\y)$, and $\y$ is unique up
to a multiplicative constant.
\end{itemize}
\end{theorem}
It follows from Theorem~\ref{thm:pf} that for any eigenpair $(\lambda, \mathbf x)$ of $\A\in N_{m,n}$ (i.e., $\A\mathbf x^{m-1}=\lambda\x^{[m-1]}$), whenever $\x\in\mathbb R^n_{++}$ we have $\lambda=\rho(\A)$.

We will call a nonnegative eigenvector of $\A\in N_{m,n}$ corresponding to $\rho(\A)$ as a \textit{nonnegative Perron vector}, and a positive eigenvector of $\A\in N_{m,n}$ as a \textit{positive Perron vector}. Thus, each nonnegative tensor has a nonnegative Perron vector, and a weakly irreducible nonnegative tensor has a unique positive Perron vector.

The next theorem will help us to prove our main theorem.
\begin{theorem}\label{thm:fixed}
Let integers $m\geq 2$ and $n\geq 2$.
Let $g_i\in\mathbb R_+[\mathbf x]$ be polynomials in $\x$ with nonnegative coefficients for all $i\in [n]$. If there are two positive
vectors $\mathbf y,\mathbf z\in \mathbb R^n_{++}$ such that
\[
\mathbf y\leq \mathbf z,\ g_i(\mathbf y)\geq y_i^{m-1}\ \text{and }g_i(\mathbf z)\leq z_i^{m-1}\ \text{for all }i\in [n],
\]
then there exists a vector $\mathbf w\in [\mathbf y,\mathbf z]:=\{\mathbf x : y_i\leq x_i\leq z_i\ \text{for all }i\in [n]\}$ such that
\[
g_i(\mathbf w)=w_i^{m-1}\ \text{for all }i\in [n].
\]
Moreover, for any initial point $\x_0\in [\y,\mathbf z]$, the iteration
\[
(\x_{k+1})_i:=\big[g_i(\x_k)\big]^{\frac{1}{m-1}}\ \text{for all }i\in [n]
\]
satisfies
\begin{enumerate}
\item $\x_{k+1}\geq \x_k$, and
\item $\lim_{k\rightarrow\infty}\x_k=\x_*$ with $\x_*\in\mathbb R^n_{++}$ such that $g_i(\x_*)=(\x_*)_i^{m-1}$ for all $i\in [n]$.
\end{enumerate}
\end{theorem}

\proof
Define $f_i : \mathbb R^n_{++}\rightarrow\mathbb R_{++}$ as
\[
f_i(\mathbf x):=\big[g_i(\mathbf x)\big]^{\frac{1}{m-1}}\ \text{for all }i\in [n].
\]
It follows from $g_i(\mathbf y)\geq y_i^{m-1}>0$ that the mapping $f:=(f_1,\dots,f_n)^\mathsf{T} : \mathbb R^n_{++}\rightarrow\mathbb R^n_{++}$ is well-defined. Since $g_i$'s are polynomials with nonnegative coefficients, the mapping $f$ is clearly increasing in the interval $[\mathbf y,\mathbf z]$, i.e., $f(\mathbf x_1)\geq f(\mathbf x_2)$ whenever $\mathbf x_1-\mathbf x_2\in\mathbb R^n_+$ and $\x_1,\x_2\in [\mathbf y,\mathbf z]$; and compact on every sub-interval of $[\mathbf y,\mathbf z]$, i.e.,  $f$ is continuous and maps sub-intervals into compact sets. Note that $f(\mathbf y)\geq \mathbf y$ and $f(\mathbf z)\leq \mathbf z$. It follows from \cite[Theorem~6.1]{A} that there exists $\mathbf w\in [\mathbf y,\mathbf z]$ such that $f(\mathbf w)=\mathbf w$, which is exactly the first half of the desired result.

With the established results, the convergence of the iteration follows also from \cite[Theorem~6.1]{A}.
\qed

\section{Nonnegative Tensor Partition}\label{sec:partition}
Given a tensor $\A\in N_{m,n}$ and an index subset $I=\{j_1,\dots,j_{|I|}\}\subseteq\{1,\dots,n\}$, $\A_I\in N_{m,|I|}$ is the $m$th order $|I|$-dimensional \textit{principal sub-tensor} of $\A$ defined as
\[
(\A_I)_{i_1\dots i_m}= a_{j_{i_1}\dots j_{i_m}}\ \text{for all }i_j\in [|I|]\ \text{and }j\in [m].
\]
In particular, $\mathbf x_I\in\mathbb R^{|I|}$ is the \textit{sub-vector} of $\mathbf x$ indexed by $I$.

Let $\mathfrak{S}(n)$ be the group of permutations on $n$ elements, or the symmetric group on the set $[n]$ in some literature.
We can define a group action on $N_{m,n}$ by $\mathfrak{S}(n)$ as:
\[
(\sigma\cdot \A)_{i_1\dots i_m}=a_{\sigma(i_1)\dots \sigma(i_m)}\ \text{for }\sigma\in\mathfrak{S}(n)\ \text{and }\A\in N_{m,n}.
\]
Let $\A\x^{m-1}=\lambda\x^{[m-1]}$. For any $\sigma\in\mathfrak{S}(n)$, define $\y\in\mathbb C^n$ with $y_i=x_{\sigma(i)}$. We have $(\sigma\cdot\A )\y^{m-1}=\lambda\y^{[m-1]}$. Therefore, for any $\A\in N_{m,n}$, tensors in the orbit $\{\sigma\cdot\A : \sigma\in\mathfrak{S}(n)\}$ have the same set of eigenvalues. In particular,
\begin{equation}\label{eqn:radius}
\rho(\sigma\cdot\A)=\rho(\A)\ \text{for all }\sigma\in\mathfrak{S}(n).
\end{equation}

The main result in \cite{HHQ} can be stated as follows.
\begin{proposition}[Nonnegative Tensor Partition]\label{prop:partition}
For any $\A\in N_{m,n}$, there exists a partition of the index set $[n]$:
\[
I_1\cup\dots\cup I_r=[n]
\]
such that for all $j=1,\dots,r$
\begin{multline}\label{eqn:partition}
\A_{I_j}\ \text{is weakly irreducible, and }
a_{si_2\dots i_m}=0\ \text{for all }s\in I_j\\ \text{and }(i_2,\dots,i_m)\in I(t)\cap \big(\cup_{k=1}^jI_k\big)^{m-1}\ \text{for all }t\in I_1\cup\dots\cup I_{j-1}.
\end{multline}
\end{proposition}
Note that each $\A_{I_j}$ is a weakly irreducible principal sub-tensor of $\A$.
We can assume that
\[
I_j=[s_j]\setminus [s_{j-1}]
\]
with $s_0=0$, $s_0<s_1<\dots<s_r=n$, and $[0]:=\emptyset$.
In general, it should be: there exists a
$\sigma\in\mathfrak{S}(n)$ such that $\sigma\cdot\A$ has such a partition. While, in view of \eqref{eqn:radius}, we may assume throughout this paper without loss of generality that $\sigma=\text{id}$, the multiplicative identity of the group $\mathfrak{S}(n)$.

The next result is proved in \cite{HHQ}.
\begin{proposition}\label{prop:radius}
Let $\A\in N_{m,n}$ be partitioned as \eqref{eqn:partition}. Then
\[
\rho(\A)=\max\{\rho(\A_{I_j}) : j\in[r]\}.
\]
\end{proposition}

\begin{definition}[Genuine Weakly Irreducible]\label{def:genuine}
A weakly irreducible principal sub-tensor $\A_{I_j}$ of $\A$ is \textit{genuine} if
\begin{equation}\label{eqn:genuine}
a_{si_2\dots i_m}=0\ \text{for all }s\in I_j\ \text{and }(i_2,\dots,i_m)\in I(t)\ \text{for all }t\in [n]\setminus I_{j}.
\end{equation}
\end{definition}
In the matrix case, genuine weakly irreducible sub-tensors are corresponding to basic
classes of \cite{Ro}.  Note that for each tensor $\A\in N_{m,n}$, it always has one genuine weakly irreducible principal sub-tensor, namely the principal sub-tensor $\A_{I_r}$ in \eqref{eqn:partition}.

With Definition~\ref{def:genuine}, we can further rearrange $I_1,\dots,I_r$ to get a partition as follows.
\begin{proposition}[Canonical Nonnegative Tensor Partition]\label{prop:genuine}
Let $\A\in N_{m,n}$. The index set $[n]$ can be partitioned into $R\cup I_{s+1}\cup\dots\cup I_r$ with $R=I_1\cup\dots \cup I_s$ such that in addition to \eqref{eqn:partition}
\begin{enumerate}
\item $\A_{I_j}$ is a genuine weakly irreducible principal sub-tensor for all $j\in\{s+1,\dots,r\}$, and
\item for each $t\in [s]$ there exist $p_t\in I_t$ and $q_t\in I_{t+1}\cup\dots\cup I_r$ such that
\[
a_{p_ti_2\dots i_m}>0\ \text{for some }(i_2,\dots,i_m)\in I(q_t).
\]
\end{enumerate}
Moreover, the partition for the genuine principal sub-tensor blocks $\A_{I_{s+1}},\dots,\A_{I_r}$ is unique up to permutation on the index sets $\{I_{s+1},\dots,I_r\}$.
\end{proposition}

\proof Suppose that we have the tensor $\A$ with a partition as in Proposition~\ref{prop:partition}.
It then follows from Definition~\ref{def:genuine} that a weakly irreducible principal sub-tensor $\A_{I_j}$ is genuine if and only if
\[
(M_{\A})_{I_j}=M_{\A_{I_j}},\ \text{i.e., }a_{i_1i_2\dots i_m}=0\ \text{whenever }i_1\in I_j\ \text{and }\{i_2\dots,i_m\}\cap I_j^\mathsf{\complement}\neq\emptyset.
\]
Therefore, the genuine weakly irreducible principal sub-tensors are uniquely determined, and
we can group the genuine weakly irreducible principal sub-tensors together, say $I_{s+1},\dots,I_r$ without loss of generality. This sorting can be done without destroying the relative back and forth orders of the blocks which are not genuine.

For any $j\in [s]$,
since $\A_{I_j}$ is not a genuine weakly irreducible principal sub-tensor, there exists a $p_j\in I_j$ and $q_j\notin I_j$ such that
\[
a_{p_ji_2'\dots i_m'}>0\ \text{for some }(i_2',\dots,i_m')\in I(q_j).
\]
However, from \eqref{eqn:partition} for all $t\in I_1\cup\dots\cup I_{j-1}$

\[
a_{si_2\dots i_m}=0\ \text{for all }s\in I_j\ \text{and }(i_2,\dots,i_m)\in I(t)\cap \big(\cup_{k=1}^jI_k\big)^{m-1},
\]
we must have that
\[
\{i_2',\dots,i_m'\}\cap (I_{j+1}\cup\dots\cup I_r)\neq\emptyset,
\]
since otherwise $q_j\in I_1\cup\dots\cup I_{j-1}$ and $(i_2',\dots,i_m')\in  \big(\cup_{k=1}^jI_j\big)^{m-1}$. Thus, $q_j$ can be chosen in $I_{j+1}\cup\dots\cup I_r$. 
The result then follows. \qed

It is easy to see that $r\geq 1$ and $s\leq r-1$, since always there is a genuine weakly irreducible principal sub-tensor. A tensor $\A\in N_{m,n}$ in the form described in Proposition~\ref{prop:genuine} is in a \textit{canonical nonnegative partition form}. It follows that, each nonnegative tensor (up to a group action by $\mathfrak{S}(n)$) can be written in a canonical nonnegative partition form.


\section{Nontrivially Nonnegative Tensors}\label{sec:nnt}
Let $\mathbf e\in\mathbb R^n$  be the vector of all ones. We will write $\mathbf x\geq \mathbf 0$ ($>\mathbf 0$ respectively) whenever $\mathbf x\in\mathbb R^n_+$ ($\mathbb R^n_{++}$ respectively).
\begin{definition}[\cite{HHQ}]\label{def:strictly}
A tensor $\A\in N_{m,n}$ is \textit{strictly nonnegative}, if $\A\mathbf e^{m-1}>\mathbf 0$.
\end{definition}
It is easy to see that in Definition~\ref{def:strictly}, $\A\in N_{m,n}$ being strictly nonnegative is equivalent to $\A\mathbf x^{m-1}>\mathbf 0$ for any positive vector $\mathbf x$.

The next observation follows from Definition~\ref{def:strictly} and Theorem~\ref{thm:pf}.
\begin{proposition}\label{prop:weak-strict}
If $\A\in N_{m,n}$ is weakly irreducible and $\rho(\A)>0$, then $\A$ is strictly nonnegative.
\end{proposition}

\begin{definition}[Nontrivially Nonnegative Tensor]\label{def:non-strictly}
A tensor $\A\in N_{m,n}$ is \textit{nontrivially nonnegative}, if there exists a nonempty subset $I\subseteq\{1,\dots,n\}$ such that $\A_I$ is strictly nonnegative.
\end{definition}

\begin{theorem}[Positive Spectral Radius]\label{thm:non-strict}
Given any $\A\in N_{m,n}$, $\rho(\A)>0$ if and only if $\A$ is nontrivially nonnegative.
\end{theorem}

\proof
Necessity: Suppose that $\mathbf x\geq \mathbf 0$ is an eigenvector of $\A$ such that (cf.\ Theorem~\ref{thm:pf})
\[
\A \mathbf x^{m-1}=\rho(\A)\mathbf x^{[m-1]}
\]
and $\rho(\A)>0$. From Section~\ref{sec:partition}, $[n]$ can be partitioned into $\{I_1,\dots,I_r\}$ such that each $\A_{I_j}$ is weakly irreducible. By Proposition~\ref{prop:radius}, $\rho(\A_{I_j})=\rho(\A)>0$ for some $j\in [r]$. Therefore, the principal sub-tensor $\A_{I_j}$ is strictly nonnegative by Proposition~\ref{prop:weak-strict}. By Definition~\ref{def:non-strictly}, $\mathcal A$ is nontrivially nonnegative.

Sufficiency: Suppose that $\A_I$ is strictly nonnegative for some index subset $I\subseteq\{1,\dots,n\}$. It follows from Definition~\ref{def:strictly} that
\[
\A_I\mathbf e_I^{m-1}>\mathbf 0,
\]
which implies that
\[
(\A\mathbf y^{m-1})_I\geq \A_I\mathbf e_I^{m-1}>\mathbf 0,
\]
where $\mathbf y_I=\mathbf e_I$ and $y_i=0$ for $i\notin I$.
By \cite[Theorem~5.5]{YY}, we have
\[
\rho(\A)=\max_{\mathbf x\geq 0}\min_{x_i>0}\frac{(\mathcal A\mathbf x^{m-1})_i}{x_i^{m-1}}\geq \min_{i\in I}(\mathcal A\mathbf y^{m-1})_i\geq \min_{i\in I}(\A_I\mathbf e_I^{m-1})_i>0.
\]
Therefore, we have $\rho(\A)>0$.
\qed

\section{Positive Perron Vector}\label{sec:positive-vector}

\begin{proposition}\label{prop:p1}
Suppose that $\A = (a_{i_1\ldots i_m}) \in N_{m, n}$ has a positive eigenvector.  Then it is a positive Perron vector and
$\A$ is either the zero tensor or a strictly nonnegative tensor.
\end{proposition}
\proof  Suppose that $\A$ has a positive eigenvector $\x$, corresponding to an eigenvalue $\lambda$. Then $\lambda=\rho(\A)$ by Theorem~\ref{thm:pf} as $\A$ is nonnegative.  If $\lambda = 0$, then by the eigenvalue equations, $\A$ must be the zero tensor.   Suppose that $\lambda > 0$.   Then by the fact that $\A\x^{m-1}>\mathbf 0$ is a positive vector, the result follows from Definition~\ref{def:strictly}.
\qed

\begin{definition}[Strongly Nonnegative Tensor]\label{def:strong}
Let $\A\in N_{m,n}$ have a canonical nonnegative partition as in Proposition~\ref{prop:genuine}. $\A$ is called \textit{strongly nonnegative}, if
\begin{equation}\label{eqn:nec-suf}
\begin{cases}\rho(A_{I_j})=\rho(\A)& \text{if }\A_{I_j}\ \text{is genuine}\\ \rho(A_{I_j})<\rho(\A)&\text{otherwise}.\end{cases}
\end{equation}
\end{definition}

With Definition~\ref{def:strong}, we present our main theorem.
\begin{theorem}[Positive Perron Vector]\label{thm:eigenvector}
Let $\A\in N_{m,n}$. Then $\A$ has a positive Perron vector if and only if $\A$ is strongly nonnegative.
\end{theorem}
Theorem~\ref{thm:eigenvector} will be proved in Section~\ref{sec:proof}, after the preparations Sections~\ref{sec:system} and \ref{sec:lemmas}.

By Theorems \ref{thm:pf}, \ref{thm:eigenvector} and Proposition \ref{prop:p1}, we see that a weakly irreducible nonnegative tensor is a strongly nonnegative tensor, and a nonzero strongly nonnegative tensor is a strictly nonnegative tensor.


\subsection{Systems of eigenvalue equations}\label{sec:system}
In the following, we will always assume that a given tensor $\A\in N_{m,n}$ is in a canonical nonnegative partition form (cf.\ Proposition~\ref{prop:genuine}). Recall that $[n]=R\cup I_{s+1}\dots\cup I_r$ with $R=I_1\cup\dots\cup I_s$.

For any $j\in [r]$, if we let $K_j:=[n] \setminus I_j$, then
\begin{equation}\label{eqn:tensor-ref}
(\A\x^{m-1})_{I_j}=\A_{I_j}\x_{I_j}^{m-1}+\sum_{u=1}^{m-1}\A_{j,u}(\x_{K_j})\x_{I_j}^{u-1}
\end{equation}
for some tensors $\A_{j,u}(\x_{K_j})\in T_{u,|I_j|}(\mathbb R_+[\x_{K_j}])$ for all $u=1,\dots,m-1$. Namely, $\A_{j,u}(\x_{K_j})$ is a tensor of order $u$ and dimension $|I_j|$ with the entries being polynomials in the variables $\x_{K_j}$ with coefficients in the set $\mathbb R_+$. Moreover, it follows from \eqref{eqn:tensor-ref} that each entry of $\A_{j,u}(\x_{K_j})$ is either zero or a homogeneous polynomial of degree $m-u$. We note that there can be many choices of
 tensors $\A_{j,u}(\x_{K_j})\in T_{u,|I_j|}(\mathbb R_+[\x_{K_j}])$ to fulfill the system \eqref{eqn:tensor-ref}, similar to the rationale that there are many tensors $\T\in T_{m,n}$ which can result in the same polynomial system $\T\x^{m-1}$. However, it is well-defined in the sense that the polynomial systems $\A_{j,u}(\x_{K_j})\x_{I_j}^{u-1}$'s are all uniquely determined by $\A$. We note that when the system of polynomials $\A_{j,u}(\x_{K_j})\x_{I_j}^{u-1}=\mathbf 0$, the tensor $\A_{j,u}(\x_{K_j})\in T_{u,|I_j|}(\mathbb R_+[\x_{K_j}])$ is uniquely determined as the zero tensor.
\begin{lemma}\label{lemma:genuine}
Suppose the notation is adopted as above.  Then a weakly irreducible principal sub-tensor $\A_{I_j}$ of $\A$ is genuine if and only if
\begin{equation}\label{eqn:nonzero}
\sum_{u=1}^{m-1}\A_{j,u}(\mathbf e_{K_j})\mathbf e_{I_j}^{u-1}=\mathbf 0
\end{equation}
is the zero vector, which is further equivalent to each tensor $\A_{j,u}(\mathbf x_{K_j})$ is the zero tensor for all $u\in [m-1]$.
\end{lemma}

\proof
It follows from \eqref{eqn:genuine} that a weakly irreducible principal sub-tensor $\A_{I_j}$ is genuine if and only if the right polynomials of $\x$ in \eqref{eqn:tensor-ref} only involve variables $\{x_t : t\in I_j\}$. This is equivalent to $\A_{j,u}(\x_{K_j})=0$ for all $u=1,\dots,m-1$ for any choice of $\A_{j,u}(\x_{K_j})$ in \eqref{eqn:tensor-ref}.

From the facts
\begin{itemize}
\item the tensors $\A_{j,u}(\x_{K_j})$ all take polynomials with nonnegative coefficients as entries, and
\item $\A_{j,u}(\x_{K_j})\x_{I_j}^{u-1}$'s are uniquely determined,
\end{itemize}
we have that the above zero polynomials condition is equivalent to each entry of all the tensors $\A_{j,u}(\x_{K_j})$ is zero. This is further equivalent to
\[
\sum_{u=1}^{m-1}\A_{j,u}(\mathbf e_{K_j})\mathbf e_{I_j}^{u-1}=\mathbf 0.
\]
\qed

For all $j\in [s-1]$, let $L_j=R\setminus I_j$.
From \eqref{eqn:partition},
we can further partition $\A_{j,u}(\x_{K_j})$ for all $j\in [s]$ into two parts
\begin{equation}\label{eqn:tensor-sub}
\A_{j,u}(\x_{K_j})=\mathcal H_{j,u}(\x_{L_j})+\B_{j,u}(\x_{K_j})\ \text{for all }u\in [m-1],
\end{equation}
where $\mathcal H_{j,u}(\x_{L_j})\in T_{u,|I_j|}(\mathbb R_+[\x_{L_j}])$ with each entry being either zero or a homogeneous polynomial of degree $m-u$ in the variables $\x_{L_j}$, and $\B_{j,u}(\x_{K_j})\in T_{u,|I_j|}(\mathbb R_+[\x_{K_j}][\x_{L_j}])$ with each entry being either zero or a polynomial of degree in the variables $\x_{L_j}$ strictly smaller than $m-u$.


\begin{proposition}\label{prop:non-genuine}
Suppose the notation is adopted as above.  Then
\begin{itemize}
\item [(a)] for $j\in [s]$ it follows $\sum_{u=1}^{m-1}\mathcal H_{j,u}(\mathbf y_{L_j})\mathbf e_{I_j}^{u-1}=\mathbf 0$ with $\y_{I_t}=\mathbf e_{I_t}$ for $t\in [j-1]$ and $\y_{I_t}=\mathbf 0$ for the others $t\in [s]\setminus\{j\}$,
\item [(b)] $\sum_{u=1}^{m-1}\B_{s,u}(\mathbf e_{K_s})\mathbf e_{I_s}^{u-1}\neq \mathbf 0$, and
\item [(c)] for $j\in [s-1]$, either $\sum_{u=1}^{m-1}\B_{s,u}(\mathbf e_{K_s})\mathbf e_{I_s}^{u-1}\neq \mathbf 0$, or
$\sum_{u=1}^{m-1}\mathcal H_{j,u}(\mathbf e_{L_j})\mathbf e_{I_j}^{u-1}\neq \mathbf 0$.
\end{itemize}
\end{proposition}

\proof
Note that item (a) follows from \eqref{eqn:partition}, and from which $\mathcal H_{s,u}(\y_{L_s})=\mathbf 0$ for all $u\in [m-1]$.
Items (b) and (c) then follow from Proposition~\ref{prop:genuine} that for each $j\in [s]$ there exists $p_j\in I_j$ and $q_j\in I_{j+1}\cup\dots\cup I_r$ such that
\[
a_{p_ji_2\dots i_m}>0\ \text{for some }(i_2,\dots,i_m)\in I(q_j),
\]
which implies that
\[
\sum_{u=1}^{m-1}\mathcal A_{j,u}(\mathbf e_{K_j})\mathbf e_{I_j}^{u-1}\neq \mathbf 0.
\]
The results now follow.
\qed

\subsection{Solvability of polynomial systems}\label{sec:lemmas}
The notation in the sequel is independent of the previous one.
\begin{lemma}\label{lemma:weak-positive}
Let $\A\in N_{m,n}$. For arbitrary $\epsilon>0$, there exists a positive vector $\x\in\mathbb R^n_{++}$ such that
\[
\A\x^{m-1}\leq(\rho(\A)+\epsilon)\x^{[m-1]}.
\]
\end{lemma}

\proof
Suppose that $\A$ is weakly reducible, since otherwise the conclusion follows from Theorem~\ref{thm:pf}.
Thus, we assume that $I_1\cup\dots \cup I_r=[n]$ is a partition of $\A$ as Proposition~\ref{prop:partition} (cf.\ \ref{eqn:partition}).

The proof is by induction on the block number $r$. The case when $r=1$ follows from Theorem~\ref{thm:pf} as we showed.  
Suppose that the conclusion is true when $r=s-1$ for some $s\geq 2$. In the next, we assume that $r=s$. Let $\kappa=\frac{1}{2}\epsilon$. 
Denote by $\mathcal C=\mathcal A_{I_1\cup\dots\cup I_{r-1}}\in N_{m,n-|I_r|}$ the principal sub-tensor of $\A$. 
It is easy to see that $I_1\cup\dots\cup I_{r-1}$ is a partition of $\mathcal C$ as Proposition~\ref{prop:partition}. Therefore, by the inductive hypothesis, we can find a vector $\y\in\mathbb R^{|I_1|+\dots+|I_{r-1}|}_{++}$ such that
\[
\mathcal C\y^{m-1}\leq(\rho(\mathcal C)+\kappa)\y^{[m-1]}\leq(\rho(\A)+\kappa)\y^{[m-1]},
\]
where the last inequality follows from $\rho(\mathcal C)=\max\{\rho(\A_{I_j}) : j\in [r-1]\}\leq \max\{\rho(\A_{I_j}) : j\in [r]\}=\rho(\A)$ by Proposition~\ref{prop:radius}. 
We also have that there exists $\mathbf z\in\mathbb R^{|I_r|}_{++}$ such that
\[
\A_{I_r}\mathbf z^{m-1}\leq(\rho(\A_{I_r})+\kappa)\mathbf z^{[m-1]}. 
\]
It follows from Proposition~\ref{prop:partition} that there are some tensors $\mathcal C_u(\mathbf z)\in T_{u,n-|I_r|}(\mathbb R_+[\mathbf z])$ for $u\in [m-1]$ such that with $\mathbf w:=(\beta\mathbf y^\mathsf{T},\mathbf z^\mathsf{T})^\mathsf{T}\in\mathbb R^n$
\[
(\A\mathbf w^{m-1})_{I_1\cup\dots\cup I_{r-1}}=\beta^{m-1}\mathcal C\y^{m-1}+\sum_{u=1}^{m-1}\beta^{u-1}\mathcal C_u(\mathbf z)\y^{u-1}\ \text{and }(\A\mathbf w^{m-1})_{I_r}=\A_{I_r}\mathbf z^{m-1}.
\]
Thus, when $\beta>0$ is sufficiently large we have
\[
(\A\mathbf w^{m-1})_{I_1\cup\dots\cup I_{r-1}}\leq \beta^{m-1}(\rho(\A)+2\kappa)\y^{[m-1]}=(\rho(\A)+2\kappa)\mathbf w_{I_1\cup\dots\cup I_{r-1}}^{[m-1]}=(\rho(\A)+\epsilon)\mathbf w_{I_1\cup\dots\cup I_{r-1}}^{[m-1]}
\]
as well as 
\[
(\A\mathbf w^{m-1})_{I_r}=\A_{I_r}\mathbf z^{m-1}\leq(\rho(\A_{I_r})+\kappa)\mathbf z^{[m-1]}\leq(\rho(\A_{I_r})+\epsilon)\mathbf w_{I_r}^{[m-1]}\leq (\rho(\A)+\epsilon)\mathbf w_{I_r}^{[m-1]}
\]
The result then follows. 
\qed

\begin{lemma}\label{lemma:suf}
Let $\lambda>0$, integers $n, s>0$ and partition $I_1\cup\dots\cup I_s=[n]$.
Suppose that for all $j\in [s]$, $\A_{I_j}\in N_{m,|I_j|}$ is weakly irreducible with $\rho(\A_{I_j})<\lambda$,
and $\A_{j,u}(\mathbf x)\in T_{u,|I_j|}(\mathbb R_+[\mathbf x_{[n]\setminus I_j}])$ for $u=1,\dots,m-1$ are such that
\begin{enumerate}
\item the degree of each entry of $\A_{j,u}(\mathbf x)$ is not greater than $m-u$, 
\item if let $\B_{j,u}(\mathbf x)\in T_{u,|I_j|}(\mathbb R_+[\mathbf x_{[n]\setminus I_j}])$ by deleting polynomials of degree $m-u$ in each entry of $\A_{j,u}(\x)$, and $\mathcal H_{j,u}(\mathbf x)=\mathcal A_{j,u}(\mathbf x)-\mathcal B_{j,u}(\mathbf x)$ for all $u\in [m-1]$ and $j\in [s]$, then $\mathcal H_{j,u}(\mathbf w^{(j)})\mathbf e_{I_j}^{u-1}=\mathbf 0$ with $\mathbf w^{(j)}_{I_1\cup\dots\cup I_{j-1}}=\mathbf e_{I_1\cup\dots\cup I_{j-1}}$ and $\mathbf w^{(j)}_{I_j\cup\dots\cup I_s}=\mathbf 0$, for all $j\in [s]$, and
\item it holds
\begin{equation}\label{eqn:lemma-connect}
\sum_{u=1}^{m-1}\B_{s,u}(\mathbf e)\mathbf e_{I_s}^{u-1}\neq \mathbf 0,
\end{equation}
and with $\y=\mathbf e_{I_1\cup\dots\cup I_j}+t\mathbf e_{I_{j+1}\cup\dots\cup I_s}$
\begin{equation}\label{eqn:lemma-nonzero}
\sum_{u=1}^{m-1}\B_{j,u}(\mathbf e)\mathbf e_{I_j}^{u-1}\neq \mathbf 0,\  \text{or }
\lim_{t\rightarrow\infty}\bigg\|\sum_{u=1}^{m-1}\A_{j,u}(\y)\mathbf e_{I_j}^{u-1}\bigg\|\to\infty
\end{equation}
for all $j\in [s-1]$.
\end{enumerate}
Then we have that there is a positive solution $\x\in\mathbb R^n_{++}$ for the following system
\begin{equation}\label{eqn:positive}
\A_{I_j}\x_{I_j}^{m-1}+\sum_{u=1}^{m-1}\A_{j,u}(\x)\x_{I_j}^{u-1}=\lambda\x_{I_j}^{[m-1]}\ \text{for all }j\in [s].
\end{equation}
\end{lemma}

\proof
We divide the proof into three parts.

{\bf Part I.}
Let $f:=(f_{I_1},\dots,f_{I_s}) : \mathbb R^n_{++}\rightarrow\mathbb R^n_{++}$ with
\[
f_{I_j}(\mathbf x):=\left[\frac{1}{\lambda}\left(\A_{I_j}\x_{I_j}^{m-1}+\sum_{u=1}^{m-1}\A_{j,u}(\x)\x_{I_j}^{u-1}\right)\right]^{\frac{1}{[m-1]}}.
\]
Since $\A_{I_j}$ is weakly irreducible and $\A_{j,u}(\mathbf x)$'s are tensors with entries being nonnegative polynomials, $f_{I_j}: \mathbb R^n_{++}\to\mathbb R^{|I_j|}_{++}$ is well-defined when either $|I_j|>1$ or $\A_{I_j}>0$ when $|I_j|=1$. The case when $|I_j|=1$ and $\A_{I_j}=0$ is also well-defined, since \eqref{eqn:lemma-nonzero} implies the existence of a positive entry. Therefore, the map $f : \mathbb R^n_{++}\to\mathbb R^n_{++}$ is well-defined.

{\bf Part II.}
Let $\mathcal A\in N_{m,n}$ be the tensor with the principal sub-tensors $\A_{I_j}$ for $j\in [s]$
and such that it satisfies the polynomial systems
\[
(\mathcal A\mathbf x^{m-1})_{I_j}=\mathcal A_{I_j}\mathbf x_{I_j}^{m-1}+\sum_{u=1}^{m-1}\mathcal H_{j,u}(\mathbf x)\mathbf x_{I_j}^{u-1}\ \text{for all }j\in [s].
\]
It follows from the second listed hypothesis and Proposition~\ref{prop:partition} that $I_1\cup\dots\cup I_s=[n]$ forms a partition for the tensor $\mathcal A$.  By Proposition~\ref{prop:radius}, $\rho(\A)=\max\{\rho(\A_{I_j}) : j\in [s]\}<\lambda$. 

Since $\lambda>\rho(\A)$, it follows from Lemma~\ref{lemma:weak-positive} that there exists a vector $\mathbf y>\mathbf 0$ such that
\begin{equation}\label{induction}
\mathcal A\mathbf y^{m-1}<\lambda\mathbf y^{[m-1]}. 
\end{equation}
So, with $\beta>0$, we have $\beta\y>\mathbf 0$ and
\begin{align}\label{upper}
f_{I_j}(\beta\mathbf y)&=\beta\left[\frac{1}{\lambda}\left(\mathcal A_{I_j}\mathbf y_{I_j}^{m-1}+\sum_{u=1}^{m-1}\beta^{u-m}\mathcal A_{j,u}(\beta\y)\mathbf y_{I_j}^{u-1}\right)\right]^{\frac{1}{[m-1]}}\nonumber\\
&=\beta\left[\frac{1}{\lambda}\left(\mathcal A_{I_j}\mathbf y_{I_j}^{m-1}+\sum_{u=1}^{m-1}\beta^{u-m}\mathcal H_{j,u}(\beta\y)\mathbf y_{I_j}^{u-1}+\sum_{u=1}^{m-1}\beta^{u-m}\mathcal B_{j,u}(\beta\y)\mathbf y_{I_j}^{u-1}\right)\right]^{\frac{1}{[m-1]}}\nonumber\\
&=\beta\left[\frac{1}{\lambda}\left((\mathcal A\y^{m-1})_{I_j}+\sum_{u=1}^{m-1}\beta^{u-m}\mathcal B_{j,u}(\beta\y)\mathbf y_{I_j}^{u-1}\right)\right]^{\frac{1}{[m-1]}}\nonumber\\
&\leq\beta\y_{I_j}
\end{align}
for sufficiently large $\beta>0$. Here, the inequality follows from \eqref{induction} and the fact that the maximal possible degree for the polynomials in the entries of each tensor $\mathcal B_{j,u}(\y)$ is $m-u-1$ for all $u\in [m-1]$. 
Since there are finite $j$'s, $f(\beta\y)\leq \beta \y$ for some sufficiently large $\beta$.

{\bf Part III.}
Recall that
$\B_{j,u}(\mathbf x)\in T_{u,|I_j|}(\mathbb R_+[\mathbf x_{[n]\setminus I_j}])$ is obtained by deleting polynomials of degree $m-u$ in each entry of $\A_{j,u}(\x)$ for all $u\in [m-1]$ and $j\in [s]$.
Let
\[
P_j:=\text{supp}\left(\sum_{u=1}^{m-1}\B_{j,u}(\mathbf e)\mathbf e_{I_j}^{i-1}\right)\subseteq I_j
\]
for all $j\in [s]$, and for $j\in [s-1]$
\begin{multline}\label{eqn:set-q}
Q_j:=\bigg\{z\in I_j : \lim_{t\rightarrow\infty}\sum_{w\in I_{j+1}\cup\dots\cup I_s}\bigg(\sum_{u=1}^{m-1}\A_{j,u}(\x_w)\mathbf e_{I_j}^{u-1}\bigg)_z\to\infty,\\ \text{with }(\x_w)_w=t\ \text{and }(\x_w)_v=1\ \text{for the others}\bigg\}.
\end{multline}
It follows from the \eqref{eqn:lemma-connect} that $P_s\neq\emptyset$, and \eqref{eqn:lemma-nonzero} that $P_j\cup Q_j\neq\emptyset$ for $j\in [s-1]$. Let $Q_s=\emptyset$. Let $W_j:=Q_j\setminus P_j$ for $j\in [s]$.


For each $j\in [s]$, let the majorization matrix for $\A_{I_j}$ be $M_j\in\mathbb R^{|I_j|\times |I_j|}_+$. It follows from the weak irreducibility that the directed graph $G_j=(V_j=I_j,E_j)$ associated to $M_j$ is strongly connected for every $j\in [s]$.
Therefore, for any nonempty proper subset $K_j\subset I_j$ and $t\in I_j\setminus K_j$, there should be a directed path from $t$ to some $w\in K_j$, and the intermediate vertices in this path all come from the set $I_j\setminus K_j$.
We will generate a forest (a union of trees) $T=(I_1\cup\dots \cup I_s, F)$ through the following procedure:
\begin{algorithm}[Forest Generating Algorithm]\label{pro}
The input is the directed graphs $G_j$ and the sets $P_j$ and $W_j$ for $j\in [s]$.
\begin{itemize}
\item [] Step 0: Set $F=\emptyset$, $j=s$.
\item []Step 1: If $j=0$, stop; otherwise set $J_j=P_j\cup W_j$. For each $v\in W_j$, pick a $w\in I_{j+1}\cup\dots\cup I_s$ such that
$\lim_{t\rightarrow\infty}\big(\sum_{u=1}^{m-1}\A_{j,u}(\x_w)\mathbf e_{I_j}^{u-1}\big)_v\to\infty$, and add $(v,w)$ into $F$.
Go to Step 2.
\item []Step 2: Let $K_j=I_j\setminus J_j$, $S_j=\emptyset$. If $K_j=\emptyset$, go to Step 4; otherwise, go to Step 3.
\item []Step 3: Pick a vertex $v\in K_j\setminus S_j$, add a directed path in $G_j$ from $v$ to some $w\in J_j$ with all intermediate vertices being distinct and in $K_j$ into $T$, add all the vertices in this path from $K_j$ into $S_j$, go to Step 4.
\item []Step 4: If $S_j=K_j$, go to Step 6; otherwise,  go to Step 5.
\item []Step 5: If there is $v\in K_j\setminus S_j$ such that $(v,w)\in E_j$ for some $w\in S_j$, put $v$ into $S_j$ and $(v,w)$ into $F$, go to Step 4; otherwise go to Step 3,
\item []Step 6: Set $j=j-1$, go to Step 1.
\end{itemize}
\end{algorithm}

We take a short break to show that the above procedure is well-defined.
\begin{itemize}
\item Step 1 is well-defined, since $W_j\subseteq Q_j$ and $Q_j$ is defined as \eqref{eqn:set-q}.
\item Note that $P_j\cup W_j=P_j\cup Q_j\neq \emptyset$ for all $j\in [s]$.
\item Step 3 is well-defined from the words before the procedure as well as Step 5.
\end{itemize}
Since $G_j$ is strongly connected for all $j\in [s]$, the above procedure should terminate in finitely many steps.
We note that the generated forest may not be unique. For every edge $(v,w)\in F$, the vertex $v$ is a \textit{child} of the vertex $w$, and which is the \textit{parent} of the vertex $v$. A vertex with no child is a \textit{leaf}, and a vertex with no parent is a \textit{root}.
An isolated vertex is both a leaf and a root.
It is easy to see from the above procedure that every root is a vertex in $\cup_{j=1}^sP_j$, and vice verse. It is also a matter of fact that from every vertex we can get a unique root along the directed edges. Therefore, we can define the \textit{height} of a vertex unambiguously as the length of the unique directed path from it to the root. Thus, a root has height $1$. The maximum height of the vertices in a tree is the \textit{height of the tree}, and the maximum height of the trees in a forest is the \textit{height of the forest}. We denote by $h(T)$ the height of the forest $T=(I_1\cup\dots \cup I_s, F)$ generated by Algorithm~\ref{pro}.

Let $\mathbf x$ be a positive vector, $\gamma>0$, and $(v,w)\in F$. Obviously, $v$ is not a root. Suppose that $v\in I_j$. If $w\in I_j$, we have
$(v,w)\in E_j$ and
\begin{align*}\label{lower31}
(f_{I_j}(\gamma\mathbf x))_v&=\gamma\Bigg[\frac{1}{\lambda}\bigg(\sum_{i_2,\dots,i_m\in I_j}a_{vi_2\dots i_m}x_{i_2}\dots x_{i_m}+\bigg(\sum_{p=1}^{m-1}\gamma^{p-m}\mathcal A_{j,p}(\gamma\x)\mathbf x_{I_j}^{p-1}\bigg)_v\bigg)\Bigg]^{\frac{1}{m-1}}\nonumber\\
&\geq\gamma \left(\frac{1}{\lambda}a_{vi_2'\dots i_m'}x_{i_2'}\dots x_{i_m'}\right)^{\frac{1}{m-1}}
\end{align*}
for some $\{i_2',\dots,i_m'\}\in I_j^{m-1}$ such that
\[
a_{vi_2'\dots i_m'}>0 \ \text{and }w\in \{i_2',\dots,i_m'\},
\]
since $(v,w)\in F\cap I_j^2\subseteq E_j$.
Therefore, when $x_w$ is sufficiently large $(f_{I_j}(\gamma\mathbf x))_v\geq \gamma x_v$.
If $w\notin I_j$,  then $w\in I_{j+1}\cup\dots\cup I_s$ is such that (cf.\ Algorithm~\ref{pro})
\[
\lim_{t\rightarrow\infty}\big(\sum_{u=1}^{m-1}\A_{j,u}(\y)\mathbf e_{I_j}^{u-1}\big)_v\to\infty\ \text{with }y_w=t\ \text{and }y_p=1\ \text{for the other } p\in [n]\setminus\{w\}.
\]
Moreover, we should have that $v\in W_j$ (cf.\ Algorithm~\ref{pro})
and therefore,
\[
\big(\sum_{u=1}^{m-1}\B_{j,u}(\mathbf e)\mathbf e_{I_j}^{u-1}\big)_v= 0.
\]
Thus,
\[
\big(\sum_{u=1}^{m-1}\A_{j,u}(\mathbf x)\mathbf x_{I_j}^{u-1}\big)_v=\big(\sum_{u=1}^{m-1}\mathcal H_{j,u}(\mathbf x)\mathbf x_{I_j}^{u-1}\big)_v
\]
is a homogeneous polynomial of degree $m-1$.
Henceforth, if $x_w$ is sufficiently large, we have
\begin{align*}\label{lower31}
(f_{I_j}(\gamma\mathbf x))_v&=\gamma\left[\frac{1}{\lambda}\left(\sum_{i_2,\dots,i_m\in I_j}a_{vi_2\dots i_m}x_{i_2}\dots x_{i_m}+\big(\sum_{p=1}^{m-1}\gamma^{p-m}\mathcal A_{j,p}(\gamma\x)\mathbf x_{I_j}^{p-1}\big)_v\right)\right]^{\frac{1}{m-1}}\nonumber\\
&=\gamma\left[\frac{1}{\lambda}\left(\sum_{i_2,\dots,i_m\in I_j}a_{vi_2\dots i_m}x_{i_2}\dots x_{i_m}+\big(\sum_{p=1}^{m-1}\mathcal H_{j,p}(\x)\mathbf x_{I_j}^{p-1}\big)_v\right)\right]^{\frac{1}{m-1}}\nonumber\\
&\geq\gamma \left(\frac{1}{\lambda}\big(\sum_{p=1}^{m-1}\mathcal H_{j,p}(\x)\mathbf x_{I_j}^{p-1}\big)_v\right)^{\frac{1}{m-1}}\nonumber\\
&\geq\gamma x_v.
\end{align*}

To get a desired vector $\x\in\mathbb R^n_{++}$ such that $f(\x)\geq \x$, we can start with $\x=\mathbf e$ and leaves with height $h(T)$. The case $h(T)=1$ is trivial.
Suppose that $h(T)\geq 2$ and $L\subset (I_1\cup\dots\cup I_s)\setminus (P_1\cup\dots\cup P_s)$ is the set of leaf vertices of height $h(T)$. Then, we can set the parents of these leaves sufficiently large such that
\begin{equation}\label{eqn:leaf}
(f(\gamma\mathbf x))_v\geq \gamma x_v\ \text{for all }v\in L.
\end{equation}
Secondly, let us consider the set $L'$ of vertices with height $h(T)-1$ if $h(T)>2$, which includes the parents of $L$. Vertices in $L'$ are not roots. If we set the set $P'$ of the parents of vertices in $L'$ sufficiently large, we can get
\[
(f(\gamma\mathbf x))_p\geq \gamma x_p\ \text{for all }p\in L'.
\]
It follows from the above analysis that we still withhold \eqref{eqn:leaf} when we increase $x_{p'}$ for $p'\in P'$ if necessary.
The next step is to consider the set $L''$ of vertices with height $h(T)-2$ if $h(T)>3$, which includes the parents of $L'$.
In this way, $(f(\gamma\mathbf x))_v\geq \gamma x_v$ for all child vertices $v\in (I_1\cup\dots\cup I_s)\setminus (P_1\cup\dots\cup P_s)$ by increasing their parents sufficiently large successively from vertices of height $h(T)$ to vertices of height $2$.
Since we have the constructed forest structure and any $v\in (I_1\cup\dots\cup I_s)\setminus (P_1\cup\dots\cup P_s)$ is a child of some parent $w\in I_1\cup\dots\cup I_s$, we can terminate the procedure in $h(T)-1$ steps, and therefore get that
\[
(f(\gamma\mathbf x))_v\geq \gamma x_v\ \text{for all }v\in  (I_1\cup\dots\cup I_s)\setminus (P_1\cup\dots\cup P_s)
\]
for some positive $\mathbf x$. Note that, we still have the freedom to choose $\gamma>0$.

If $w\in P_j\subset P_1\cup\dots\cup P_s$ is a root, then
\begin{equation}\label{eqn:nonzero-pf}
\bigg(\sum_{u=1}^{m-1}\B_{j,u}(\mathbf e)\mathbf e_{I_j}^{i-1}\bigg)_w>0
\end{equation}
by definition. We have
\begin{align*}\label{lower31}
(f_{I_j}(\gamma\mathbf x))_w&=\gamma\left[\frac{1}{\lambda}\left(\sum_{i_2,\dots,i_m\in I_j}a_{wi_2\dots i_m}x_{i_2}\dots x_{i_m}+\big(\sum_{p=1}^{m-1}\gamma^{p-m}\mathcal A_{j,p}(\gamma\x)\mathbf x_{I_j}^{p-1}\big)_w\right)\right]^{\frac{1}{m-1}}\nonumber\\
&=\gamma\Bigg[\frac{1}{\lambda}\Big(\sum_{i_2,\dots,i_m\in I_j}a_{wi_2\dots i_m}x_{i_2}\dots x_{i_m}\nonumber\\
&\ \ \ \ \ +\big(\sum_{p=1}^{m-1}\mathcal H_{j,p}(\mathbf x)\mathbf x_{I_j}^{p-1}+\sum_{p=1}^{m-1}\gamma^{p-m}\B_{j,p}(\gamma\mathbf x)\mathbf x_{I_j}^{p-1}\big)_w\Big)\Bigg]^{\frac{1}{m-1}}\nonumber\\
&\geq\gamma \Big(\frac{1}{\lambda}\big(\sum_{p=1}^{m-1}\gamma^{p-m}\B_{j,p}(\gamma\mathbf x)\mathbf x_{I_j}^{p-1}\big)_w\Big)^{\frac{1}{m-1}}.
\end{align*}
Note that the highest degree of entries in $\B_{j,p}(\mathbf x)$ is smaller than $m-p-1$ for all $p\in [m-1]$. This, together with \eqref{eqn:nonzero-pf}, implies that
the leading term of
\[
\big(\sum_{p=1}^{m-1}\gamma^{p-m}\B_{j,p}(\gamma\mathbf x)\mathbf x_{I_j}^{p-1}\big)_w
\]
is a term of $\frac{1}{\gamma^u}$ with positive coefficient for some integer $u>0$. Therefore, if $\gamma>0$ is sufficiently small, we have
\[
(f_{I_j}(\gamma\mathbf x))_w\geq \gamma x_w
\]
for sure. Since there are only finitely many roots, we have
\[
(f(\gamma\x))_w\geq \gamma x_w\ \text{for all }w\in P_1\cup\dots\cup P_s.
\]
Therefore, we can find a $\mathbf x$ with $\gamma>0$ such that $f(\gamma\mathbf x)\geq \gamma\mathbf x$ and $\gamma\x\leq\beta\y$ (cf.\ $\beta\y$ from \textbf{Part II.}).

In summary,
\[
f(\gamma\mathbf x)\geq \gamma\mathbf x\ \text{and } f(\beta\y)\leq \beta\y.
\]
It then follows from Theorem~\ref{thm:fixed} that there is a positive $\mathbf w\in[\gamma \mathbf x,\beta\mathbf y]$ such that
\[
f(\mathbf w)=\mathbf w.
\]
It is nothing but a positive solution $\mathbf w $ to \eqref{eqn:positive}.
\qed


\begin{lemma}\label{lemma:nec}
Suppose that $\A\in N_{m,n}$ is weakly irreducible, and $\A_i\in N_{i,n}$ for $i=1,\dots,m-1$ are such that
\[
\sum_{i=1}^{m-1}\A_i\mathbf e^{i-1}\neq \mathbf 0.
\]
If for some $\lambda>0$, there is a positive solution $\x\in\mathbb R^n_{++}$ for the following system
\[
\A\mathbf x^{m-1}+\sum_{i=1}^{m-1}\A_i\mathbf x^{i-1}=\lambda\mathbf x^{[m-1]},
\]
then $\rho(\A)<\lambda$.
\end{lemma}

\proof
Suppose, without loss of generality, that
\[
I:=\{1,\dots,r\}:=\text{supp}\left(\sum_{i=1}^{m-1}\A_i\mathbf e^{i-1}\right)
\]
for some $r\leq n$. It follows from the hypothesis that
\begin{equation}\label{small}
(\A\mathbf x^{m-1})_i<\lambda x_i^{m-1}\ \text{for all }i=1,\dots,r.
\end{equation}
If $r=n$, then the result follows from Theorem~\ref{thm:pf} directly.

In the following, we assume that $r<n$. Note that $r>0$ by the assumption on $\A_i$'s. We have
\[
(\A\mathbf x^{m-1})_j=\lambda x_j^{m-1}\ \text{for all }j=r+1,\dots,n.
\]
By the weak irreducibility, there should have a $j\in J:=\{1,\dots,n\}\setminus I$ and an $i\in I$ such that
\[
a_{ji_2\dots i_m}>0\ \text{for some multiset }\{i_2,\dots,i_m\}\ni i.
\]
Therefore, there is a nonzero term in $\sum_{i_2,\dots,i_m=1}^na_{ji_2\dots i_m}x_{i_2}\dots x_{i_m}$ involving the variable $x_i$.
We can define a new positive vector, denoted also by $\mathbf x$, through decreasing $x_i$ a little bit. It follows from the nonnegativity of $\A$ that
\begin{equation}\label{new-small}
(\A\mathbf x^{m-1})_j<\lambda x_j^{m-1}.
\end{equation}
By the continuity, we can still withhold \eqref{small} for a sufficiently small decreasing of $x_i$, as well as getting \eqref{new-small}.
While, as we can see, we get at least $r+1$ strict inequalities now. Inductively in this way, we can find a positive vector $\mathbf x$ such that
\[
\mathcal A\mathbf x^{m-1}<\lambda\mathbf x^{[m-1]}.
\]
The result then follows from Theorem~\ref{thm:pf}.
\qed

\subsection{Proof of Theorem~\ref{thm:eigenvector}}\label{sec:proof}
\proof
In the proof, we assume all the notation in Section~\ref{sec:system}.
We prove the sufficiency first.

For any $j=s+1,\dots,r$, we see that
\begin{equation}\label{eqn:genuine-eig}
(\A\x^{m-1})_{I_j}=\A_{I_j}\x_{I_j}^{m-1}\ \text{for all }\x\in\mathbb C^n.
\end{equation}
It follows from Theorem~\ref{thm:pf} that there exists a positive vector $\y_j\in\mathbb R^{|I_j|}_{++}$ such that
$\A_{I_j}\y_j^{m-1}=\rho(\A_{I_j})\y_j^{[m-1]}=\rho(\A)\y_j^{[m-1]}$ for all $j=s+1,\dots,r$.

Let $\x$ be an $n$-dimensional vector with $\x_{I_j}=\y_j$ for $j=s+1,\dots,r$ and $\x_{I_1\cup\dots\cup I_s}$ indeterminantes to be determined. It is sufficient to show that the following system of polynomials has a positive solution in $\mathbb R^{|I_1|+\dots+|I_s|}_{++}$:
\begin{equation}\label{eqn:tensor-ref-2}
\A_{I_j}\x_{I_j}^{m-1}+\sum_{u=1}^{m-1}\A_{j,u}(\x_{K_j})\x_{I_j}^{u-1}=\rho(\A)\x_{I_j}^{[m-1]}\ \text{for all }j\in [s],
\end{equation}
where $\A_{j,u}(\x_{K_j})$ has the partition \eqref{eqn:tensor-sub}.

Note that $\x_{I_j}$'s are given positive vectors for $j=s+1,\dots,r$. The indeterminant variables are $\x_{I_j}$ for $j\in [s]$. Therefore, the tensor formed by the polynomials of degree $m-u$ in $\A_{j,u}(\x_{K_j})$ is $\mathcal H_{j,u}(\x_{L_j})$ for any $u\in [m-1]$ and $j\in [s]$. If for some $j\in [m-1]$, we have $\sum_{u=1}^{m-1}\B_{j,u}(\mathbf e_{K_j})\mathbf e_{I_j}^{u-1}=\mathbf 0$,  then $\sum_{u=1}^{m-1}\mathcal H_{j,u}(\mathbf e_{L_j})\mathbf e_{I_j}^{u-1}\neq \mathbf 0$ by item (c) in Proposition~\ref{prop:non-genuine}. It follows from (a) in Proposition~\ref{prop:non-genuine} that
$\sum_{u=1}^{m-1}\mathcal H_{j,u}(\mathbf h_{L_j})\mathbf e_{I_j}^{u-1}=\mathbf 0$ with $\mathbf h_{I_t}=\mathbf e_{I_t}$ for $t\in [j-1]$ and $\mathbf h_{I_t}=\mathbf 0$ for the others. Therefore, a nonzero term involving variables from $I_{j+1}\cup\dots\cup I_s$ occurs in some entry of one tensor $\mathcal H_{j,u}(\x_{L_j})$ for some $u\in [m-1]$.

In summary,
it follows from Proposition~\ref{prop:non-genuine} that \eqref{eqn:lemma-connect} and \eqref{eqn:lemma-nonzero}, as well as the second hypothesis in Lemma~\ref{lemma:suf}, are fulfilled for the system \eqref{eqn:tensor-ref-2}.
Therefore,
by Lemma~\ref{lemma:suf}, we can find a solution $\mathbf z\in\mathbb R^{|I_1|+\dots+|I_s|}_{++}$ for \eqref{eqn:tensor-ref-2}.
Therefore, $\x$ with $\x_{I_j}=\mathbf z_{I_j}$ for $j\in [s]$ and $\x_{I_j}=\y_j$ for $j=s+1,\dots,r$ is a positive Perron vector of $\A$.

For the necessity, suppose that $\A\x^{m-1}=\lambda\x^{[m-1]}$ with $\x\in\mathbb R^n_{++}$ being a positive eigenvector for some $\lambda\geq 0$. By Theorem~\ref{thm:pf}, $\lambda=\rho(\A)$. The case for $\rho(\A)=0$ is trivial. In the following, we assume $\rho(\A)>0$. We have that for each $I_j$,
\[
(\A\x^{m-1})_{I_j}=\rho(\A)\x_{I_j}^{[m-1]}.
\]
If $\A_{I_j}$ is a genuine weakly irreducible principal sub-tensor of $\A$, it follows from Proposition~\ref{prop:genuine} that
\[
\A_{I_j}\x_{I_j}^{m-1}=\rho(\A)\x_{I_j}^{[m-1]}
\]
and Theorem~\ref{thm:pf} that $\rho(\A_{I_j})=\rho(\A)$.

If $\A_{I_j}$ is not a genuine weakly irreducible principal sub-tensor of $\A$,
it follows from Proposition~\ref{prop:genuine} that
\[
\A_{I_j}\x_{I_j}^{m-1}+\sum_{u=1}^{m-1}\A_{j,u}(\mathbf x_{K_j})\mathbf x_{I_j}^{u-1}=\rho(\A)\x_{I_j}^{[m-1]}
\]
and $\sum_{u=1}^{m-1}\A_{j,u}(\mathbf e_{K_j})\mathbf e_{I_j}^{u-1}\neq\mathbf 0$;
and Lemma~\ref{lemma:nec} that $\rho(\A_{I_j})<\rho(\A)$. 
The proof is thus complete.
\qed

\section{Algorithmic Aspects}\label{sec:algoas}

In order to get a canonical nonnegative partition for a nonnegative tensor as in Proposition~\ref{prop:genuine}, we have to recursively partition the majorization matrix of the nonnegative tensor and marjorization matrices of some induced principal sub-tensors (cf.\ \cite{HHQ}).
\subsection{Majorization matrix partition}\label{sec:major}
\begin{lemma}\label{lemma:marjorization}
Let $M_{\A}=(m_{ij})$ be the majorization matrix of $\A\in N_{m,n}$. If $m_{ij}=0$ for all $i\in I$ and $j\in I^\mathsf{\complement}$ for some nonempty proper subset $I\subset [n]$, then
\begin{equation}\label{eqn:submatrix}
M_{\A_I}=\big(M_{\A}\big)_I.
\end{equation}
\end{lemma}

\proof
Without loss of generality, we can assume that $I^\mathsf{\complement}=[p]$ for some positive $p<n$. It follows from $m_{ij}=0$ for all $i>p$ and $j\in [p]$ that
\[
a_{ii_2\dots i_m}=0\ \text{for all }i>p\ \text{and }(i_2,\dots,i_m)\in I(1)\cup\dots\cup I(p).
\]
Note that for any $i'\in I$
\[
m_{ii'} = \sum_{(i_2,\dots,i_m)\in I(i')}a_{ii_2\dots i_m}= \sum_{\begin{subarray}{c}(i_2,\dots,i_m)\in I(i')\\ \{i_2,\dots,i_m\}\cap [p]=\emptyset\end{subarray}}a_{ii_2\dots i_m}
\]
where the right most summation only involves indices $\{i_2,\dots,i_m\}\subseteq I$. With the definition for the majorization matrices for nonnegative tensors, we immediately get \eqref{eqn:submatrix}.
\qed

In general, we do not simultaneously have both $M_{\A_I}=\big(M_{\A}\big)_I$ and $M_{\A_{I^\mathsf{\complement}}}=\big(M_{\A}\big)_{I^\mathsf{\complement}}$.
\begin{example}\label{exm:majorization}
Let $\A\in N_{3,3}$ with entries
\[
a_{123}=a_{213}=a_{333}=1\ \text{and }a_{ijk}=0\ \text{for the others}.
\]
Let $I=\{3\}=\{1,2\}^\mathsf{\complement}$, we have
\[
M_{\A}=\begin{bmatrix}0&1&1\\ 1&0&1\\ 0&0&1\end{bmatrix},\ M_{\A_{\{3\}}}=[1]=\big(M_{\A}\big)_{\{3\}},\ M_{\A_{\{1,2\}}}=\begin{bmatrix}0&0\\ 0&0\end{bmatrix}\neq \big(M_{\A}\big)_{\{1,2\}}=\begin{bmatrix}0&1\\1&0\end{bmatrix}.
\]
\end{example}

A partition $I_1,\dots,I_r$ of $[n]$ is a \textit{refined partition} of a partition $J_1,\dots,J_s$ of $[n]$ if
\[
J_j=I_{j_1}\cup\dots\cup I_{j_i}\ \text{for some }j_1,\dots,j_i\in [r]\ \text{for all }j\in [s].
\]

\begin{corollary}\label{corollary:algorithm2}
Let $M_{\A}=(m_{ij})$ be the majorization matrix of $\A\in N_{m,n}$. If $\A$ has a canonical nonnegative partition $\{I_1,\dots,I_r\}$, then $M_{\A}$ has an upper triangular block structure with diagonal blocks being $\{J_1,\dots,J_s\}$ for which $
\{I_1,\dots,I_r\}$ is a refined partition.
\end{corollary}

With Corollary~\ref{corollary:algorithm2}, we can get a canonical nonnegative partition for a $\A\in N_{m,n}$ by first partition its majorization matrix $M_{\A}$ into (up to permutation) an upper triangular block form with each diagonal block being irreducible, and then recursively perform the partition to each principal sub-tensor induced by these irreducible diagonal blocks. This improves the partition method proposed in \cite{HHQ}.

\subsection{Algorithms}\label{sec:algo}
Given a nonnegative tensor $\A\in N_{m,n}$, one way to find a canonical nonnegative partition as in Proposition~\ref{prop:genuine} is by recursively partition the majorization matrices of the induced principal sub-tensors (cf.\ Corollary~\ref{corollary:algorithm2}).  We will denote by \textit{Algorithm P} an algorithm which is able to find a canonical nonnegative partition for any given nonnegative tensor. There is such an algorithm (cf.\ \cite[Section~6]{HHQ}), which can be improved with Section~\ref{sec:major}.
Therefore, in the following, we will assume that we have already computed a partition $I_1\cup\dots\cup I_r=[n]$ with properties described in Proposition~\ref{prop:genuine}. Let $I_{s+1},\dots,I_r$ be the genuine weakly irreducible blocks, and $R=I_1\cup\dots\cup I_s$.

If $\A$ is weakly irreducible, we can use the following algorithm to find the spectral radius, together with the positive Perron vector in $1$-norm being one.
\begin{algorithm}[A Higher Order Power Method \cite{HHQ}]\label{algo}
The input is a strictly nonnegative tensor $\A\in N_{m,n}$.
\begin{description}
\item [Step 0] Initialization: choose $\x^{(0)}\in \mathbb R^n_{++}$.  Let $k:=0$.

\item [Step 1] Compute
\begin{eqnarray*}
\begin{array}{c}
\bar \x^{(k+1)}:=\A(\x^{(k)})^{m-1},\quad \x^{(k+1)}:=\frac{\left(\bar
\x^{(k+1)}\right)^{[\frac{1}{m-1}]}}{\mathbf e^T\left[\left(\bar
\x^{(k+1)}\right)^{[\frac{1}{m-1}]}\right]},\\
\alpha\left(\x^{(k+1)}\right):=\max_{1\leq i\leq
n}\frac{\left(\A(\x^{(k)})^{m-1}\right)_i}{\left(\x^{(k)}\right)_i^{m-1}}\quad
\mbox{\rm and}\quad \beta\left(\x^{(k+1)}\right):=\min_{1\leq i\leq
n}\frac{\left(\A(\x^{(k)})^{m-1}\right)_i}{\left(\x^{(k)}\right)_i^{m-1}}.
\end{array}
\end{eqnarray*}

\item [Step 2] If $\alpha\left(\x^{(k+1)}\right)=\beta\left(\x^{(k+1)}\right)$ or a tolerance for $\alpha\left(\x^{(k+1)}\right)-\beta\left(\x^{(k+1)}\right)$
is satisfied, stop. Otherwise, let $k:=k+1$, go to Step 1.
\end{description}
\end{algorithm}

Note that when $\A\in N_{m,n}$ is weakly irreducible, then $\A+\mathcal I$ is weakly primitive (cf.\ \cite{HHQ,ZQX}), where $\mathcal I\in N_{m,n}$ is the identity tensor. It follows that $\rho(\A)+1=\rho(\A+\mathcal I)$.
\begin{proposition}\cite[Theorem~4.1(iv)]{HHQ}\label{prop:algo-converge}
If $\B\in N_{m,n}$ is weakly irreducible, then the sequence $\{\x^{(k)}\}$ generated by Algorithm~\ref{algo} with $\A:=\B+\mathcal I$ converges globally $R$-linearly to the positive Perron vector $\x^*$ of $\B$ corresponding to $\rho(\B)$ satisfying $\mathbf e^\mathsf{T}\x^*=1$.
\end{proposition}

In the next, we present an algorithm for determining whether a nonnegative tensor $\A\in N_{m,n}$ is strongly nonnegative or not, and finding a positive Perron vector it is.
\begin{algorithm}[Positive Perron Vector Algorithm]\label{algo:total}
The input is a nonnegative tensor $\A\in N_{m,n}$.
\begin{description}
\item [Step 0] Let $\gamma$ be a given small positive scalar.
Find a canonical nonnegative partition of $\A$ by Algorithm P with genuine weakly irreducible blocks $I_{s+1},\dots,I_r$ and $R:=I_1\cup\dots\cup I_s$.

\item [Step 1] For each $j=1,\dots,r$, find the positive eigenvector $\x_j\in\mathbb R^{|I_j|}_{++}$ such that  $\A_{I_j}\x_j^{m-1}=\rho(\A_{I_j})\x_j^{[m-1]}$ by Algorithm~\ref{algo}. If $|I_j|=1$ and $\A_{I_j}=0$, then simply set $\rho(\A_{I_j})=0$ and $\x_j=1$.

\item [Step 2] If $\max\{\rho(\A_{I_j}) : j=s+1,\dots,r\}\neq \min\{\rho(\A_{I_j}) : j=s+1,\dots,r\}$, we claim that $\A$ is not strongly nonnegative and no positive Perron vector exists; stop. Otherwise, let $\lambda:=\max\{\rho(\A_{I_j}) : j=s+1,\dots,r\}$.

\item [Step 3] If $\max\{\rho(\A_{I_j}) : j\in [s]\}\geq \lambda$, we claim that $\A$ is not strongly nonnegative and no positive Perron vector exists; stop. Otherwise, let $\y\in\mathbb R^{n}_{++}$ with
\[
\y_{I_j}=\x_j\ \text{for all }j=s+1,\dots,r.
\]

\item [Step 4] Let $\mathbf w_0:=\gamma(\x_1^\mathsf{T},\dots,\x_s^\mathsf{T})^\mathsf{T}$, $\mathbf z_0:=(\mathbf w_0^\mathsf{T},\y_{R^\mathsf{\complement}}^\mathsf{T})^\mathsf{T}$ and $k=1$.

\item [Step 5] Compute
\[
\mathbf w_k:=\bigg(\frac{(\A\mathbf z_{k-1}^{m-1})_R}{\lambda}\bigg)^{\frac{1}{[m-1]}}
\]
and $\mathbf z_k:=(\mathbf w_k^\mathsf{T},\y_{R^\mathsf{\complement}}^\mathsf{T})^\mathsf{T}$.

\item [Step 6] If $\mathbf w_k=\mathbf w_{k-1}$ or a tolerance for $\|\mathbf w_k-\mathbf w_{k-1}\|_2$ is satisfied, stop. Otherwise, let $k:=k+1$, go to Step 5.
\end{description}
\end{algorithm}

It is easy to see that $\lambda=\rho(\A)$ by Proposition~\ref{prop:radius}.
\begin{proposition}\label{prop:convergence-total}
For any given $\A\in N_{m,n}$, if $\gamma$ is sufficiently small, then Algorithm~\ref{algo:total}
\begin{enumerate}
\item either terminates in Step 2 or Step 3, which concludes that $\A$ is not strongly nonnegative and there does not exist a positive Perron vector for $\A$,
\item or generates a sequence $\{\mathbf z_k\}$ such that $\mathbf z_{k+1}\geq \mathbf z_k$ and $\lim_{k\rightarrow\infty}\mathbf z_k=\mathbf z_*$ with $\mathbf z_*$ being a positive Perron vector of $\A$ (i.e., $\A\mathbf z_*^{m-1}=\rho(\A)\mathbf z_*^{[m-1]}$).
\end{enumerate}
\end{proposition}

\proof
The conclusion for termination in either Step 2 or Step 3 follows from Theorem~\ref{thm:eigenvector}.

Suppose $\rho(\A_{I_j})<\rho(\A)$ for all $j\in [s]$ and $\rho(\A_{I_j})=\rho(\A)$ for $j=s+1,\dots,r$. Then Algorithm~\ref{algo:total} will execute Steps 4-6. It follows from the proof of Lemma~\ref{lemma:suf} that we can find positive vectors $\x$ and $\mathbf z$, and sufficiently small $\kappa$ and sufficiently large $\beta$ such that
\[
\bigg(\frac{(\A\mathbf u^{m-1})_R}{\lambda}\bigg)^{\frac{1}{[m-1]}}\geq \kappa\x \ \text{and }\bigg(\frac{(\A\mathbf v^{m-1})_R}{\lambda}\bigg)^{\frac{1}{[m-1]}}\leq \beta\mathbf z,
\]
where $\mathbf u:=(\kappa\mathbf x^\mathsf{T},\y_{R^\mathsf{\complement}}^\mathsf{T})^\mathsf{T}$, $\mathbf v=(\beta\mathbf z^\mathsf{T},\y_{R^\mathsf{\complement}}^\mathsf{T})^\mathsf{T}$ and $\y$ is defined as in Step 3.  It also follows from the proof of Lemma~\ref{lemma:suf} that if $\kappa$ fulfill the above inequality, then each positive $\tau\leq \kappa$ also does. Thus,
for every sufficiently small $\gamma>0$, we can choose $\kappa$ to make sure that $\mathbf w_0\in [\kappa\x,\beta\mathbf z]$. The convergent result then follows from Theorem~\ref{thm:fixed}.
\qed

\subsection{Numerical computation}\label{sec:num}
All the experiments were conducted by MatLab on a laptop with 2.5 GHZ Intel processor with memory 4GB.
We only tested third order strongly nonnegative tensors for Algorithm~\ref{algo:total}. Details of the problems we tested are given in Examples~\ref{exm:num1}, \ref{exm:num2} and \ref{exm:num3}.
Algorithm~\ref{algo:total} is terminated whenever
\begin{equation}\label{eqn:criteria}
\|\A\mathbf z_k^2-\lambda\mathbf z_k^{[2]}\|_2<10^{-6}\ \text{and }\|\mathbf w_k-\mathbf w_{k-1}\|_2<10^{-6}.
\end{equation}
Algorithm~\ref{algo:total} successfully computed out a positive vector to satisfy the criteria \eqref{eqn:criteria} for every tested case.
Algorithm~\ref{algo} was terminated when $\alpha(\x^{(k)})-\beta(\x^{(k)})<10^{-6}$.
\begin{example}\label{exm:num1}
The tensor $\A\in N_{8,2}$ with a partition
\[
I_1:=\{1,2\},\ I_2:=\{3,4\},\ I_3:=\{5,6\},\ I_4:=\{7,8\}
\]
and a unique genuine block
$\A_4$. Note that we will simplify $\A_{I_j}$ as $\A_j$ for $j\in [4]$. The weakly irreducible principal sub-tensors are as follows.
\[
\A_1(:,:,1) =\begin{bmatrix}    0.4423 &   0.3309\\
    0.0196  &  0.4243\end{bmatrix}\ \text{and } \A_1(:,:,2) =\begin{bmatrix}
     0.2703  &  0.8217\\
    0.1971   & 0.4299\end{bmatrix}
\]
\[
\A_2(:,:,1) =\begin{bmatrix}   0.3185&    0.0900\\
    0.5341   & 0.1117\end{bmatrix}\ \text{and } \A_2(:,:,2) =\begin{bmatrix}
     0.1363 &   0.4952\\
    0.6787    &0.1897\end{bmatrix}
\]
\[
\A_3(:,:,1) =\begin{bmatrix}       0.6664 &   0.6260\\
    0.0835 &   0.6609\end{bmatrix}\ \text{and } \A_3(:,:,2) =\begin{bmatrix}
  0.7298   & 0.9823\\
    0.8908 &   0.7690\end{bmatrix}
\]
\[
\A_4(:,:,1) =\begin{bmatrix}  0.3642  &  1.0317\\
    0.6636  &  0.5388\end{bmatrix}\ \text{and } \A_4(:,:,2) =\begin{bmatrix}
    1.1045 &   1.0251\\
    0.5921 &   1.0561\end{bmatrix}
\]
The other nonzero components of the tensor $\A$ are expressed by the following equations:
\begin{align*}
(\A\mathbf x)_{I_1}&=\A_1\x_{I_1}^2+\bigg(\begin{matrix}0.8085x_1x_7\\ 0.7551x_2x_6\end{matrix}\bigg)\\
&\ \ \
+
\bigg(\begin{matrix}0.5880x_6x_7+0.1548x_3x_8 +0.1999x_4x_8  +  0.4070x_8^2\\ 0.7487x_3x_4+0.8256x_3x_6+x_4x_8\end{matrix}\bigg)\\
(\A\mathbf x)_{I_2}&=\A_2\x_{I_2}^2+\bigg(\begin{matrix} 0.5606x_3x_5+0.9296x_4x_7\\ 0\end{matrix}\bigg)\\
&\ \ \
+
\bigg(\begin{matrix}0.9009x_2x_5\\ 0.5747x_5x_8+0.8452 x_1x_6 + 0.7386x_6x_7+0.5860 x_7^2+0.2467x_8^2\end{matrix}\bigg)\\
(\A\mathbf x)_{I_3}&=\A_3\x_{I_3}^2+\bigg(\begin{matrix}0\\ 0.5801x_6x_8\end{matrix}\bigg)+
\bigg(\begin{matrix}0.1465x_3x_8+0.1891x_4x_7\\ 0.2819x_4x_7\end{matrix}\bigg)
\end{align*}
The spectral radii of $\A_1$, $\A_2$, $\A_3$ and $\A_4$ are respectively (via Algorithm~\ref{algo})
\[
\rho(\A_1)=    1.3183,\   \rho(\A_2)= 1.2581,\  \rho(\A_3)=  2.6317,\ \text{and } \rho(\A_4)=3.1253.
\]
Therefore, $\rho(\A)=\rho(\A_4)=3.1253$ by Proposition~\ref{prop:radius}. By Algorithm~\ref{algo},
we first computed out the positive Perron vectors of $\A_i$ with $1$-norm being one as $\x_i$, then we use $\mathbf w_0:=0.5(\x_1^\mathsf{T},\x_2^\mathsf{T},\x_3^\mathsf{T})^\mathsf{T}$ (i.e., $\gamma=0.5$ in Algorithm~\ref{algo:total}) as the initial point for our fixed point iteration (i.e., Step 5 and Step 6 in Algorithm~\ref{algo:total}).
It took $0.65$ seconds with $52$ iterations to compute out a positive Perron vector of $\A$
\[
\x=( 0.4462,    0.4143,    0.3808,    0.4446,    0.2943,    0.3055,   0.5257,    0.4743)^\mathsf{T}
\]
with the corresponding eigenvalue being $\rho(\A)=3.1253$.
The final residue of the eigenvalue equations is $\|\A\x^2-3.1253\x^{[2]}\|_2=9.2323\times 10^{-7}$. Figure~1 gives information on the residue of the eigenvalue equations and the distance between two successive iterations $\mathbf w_k$ and $\mathbf w_{k-1}$ in the iteration process. The magnitude is in logarithmic order.
It shows that the convergence is sublinear, as expected for the fixed point iteration.
\begin{figure}[htbp]
\centering
\includegraphics[width=6.8in]{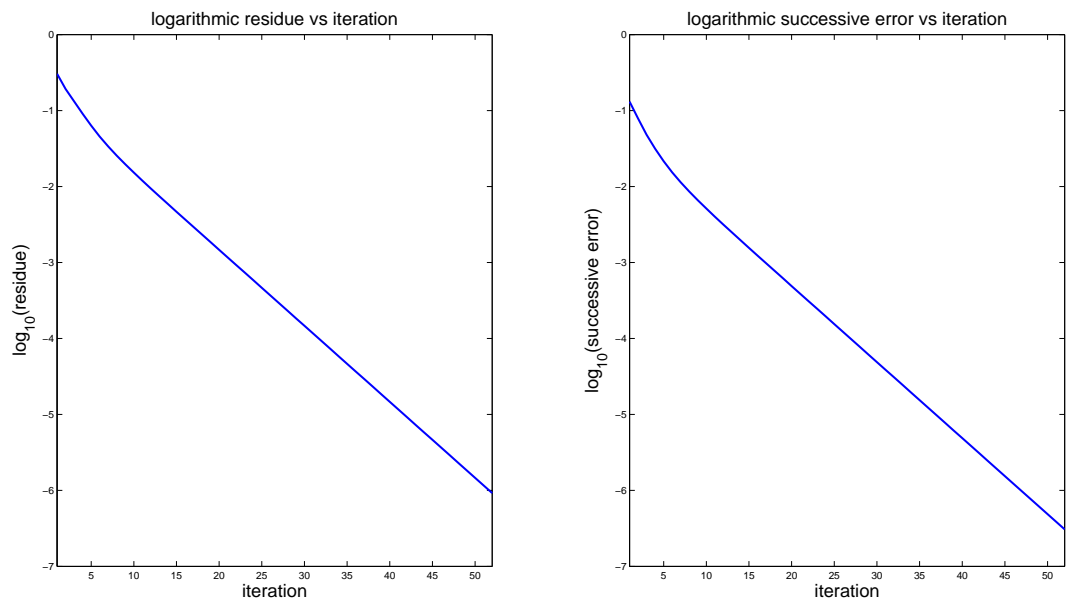}
\caption{The logarithmic residue of the eigenvalues equations $\|\A\x_k^2-\rho(\A)\x_k^{[2]}\|_2$ (left), and
the logarithmic successive error $\|\mathbf w_k-\mathbf w_{k-1}\|_2$ (right) along with the iteration.}
\end{figure}
\end{example}

\begin{example}\label{exm:num2}
In this example, we tested randomly generated nonnegative tensors. Each tested tensor has the following canonical nonnegative partition with
\[
I_1\cup I_2\cup I_3\cup I_4
\]
with a unique genuine weakly irreducible block $I_4$. The cardinality of $|I_4|$ is $10$. We tested two sets of $I_1$, $I_2$ and $I_3$.
\begin{itemize}
\item [] Case I: $|I_1|=8, |I_2|=9, |I_3|=10$, and
\item []Case II:  $|I_1|=30, |I_2|=30, |I_3|=30$.
\end{itemize}
The tensors were generated as follows
\begin{enumerate}
\item Randomly generate $\A_i\in N_{3,|I_i|}$ for $i\in [4]$. The generated tensors are all have positive components, and hence weakly irreducible.
\item Using Algorithm~\ref{algo} to compute out the spectral radii of the above generated tensors, denote the maximum of them by $\lambda_0$.
\item Let $\text{Rt}>1$ be a parameter and define $\lambda:=\lambda_0\text{Rt}$.
\item Generate the remaining components, satisfying conditions in Proposition~\ref{prop:genuine}, of each system of equations $(\A\mathbf x^2)_{I_j}$ for $j\in [3]$ randomly with sparsity $\text{den}=10\%$.
\end{enumerate}
Then, we implemented Algorithm~\ref{algo:total} to compute a positive Perron vector of the generated tensor $\A$ with the new $\A_{I_4}$ being $\frac{\lambda}{\rho(\A_4)}\A_4$. Thus, $\A$ is strongly nonnegative, and the hypothesis in Theorem~\ref{thm:eigenvector} is satisfied.
We tested various $\text{Rt}$ according to the following rule:
\begin{itemize}
\item [] Case I: $\text{Rt}(i):=1.1+0.2i$ for all $i\in [50]$, and
\item [] Case II: $\text{Rt}(i):=2+0.2i$ for all $i\in [50]$.
\end{itemize}
We will refer to $\text{Rt}$ as a measure of the ratio of the spectral radii.
For both the cases, the parameter $\gamma$ in Algorithm~\ref{algo:total} was chosen as
\[
\gamma:=10^{-5}(1/\text{Rt}(i))  \ \text{for }i\in [50].
\]
For each case and each $i\in [50]$, we made ten simulations, and recorded the means and standard variances of the numbers of iterations and the CPU time. Figure~2 contains the information for Case I, whereas Figure~3 for Case II. In each figure, we have a sub-window for an enlarged part of the whole curve. We see from the figures that it is more efficient to compute a positive Perron vector $\x$ for larger ratio of the spectral radii. This is reasonable, with larger $\text{Rt}$, $\max\{|x_i| : i\in I_1\cup I_2\cup I_3\}$ is smaller. Thus, $\x$ is closer to the initial point.
\begin{figure}[htbp]
\centering
\includegraphics[width=6.8in]{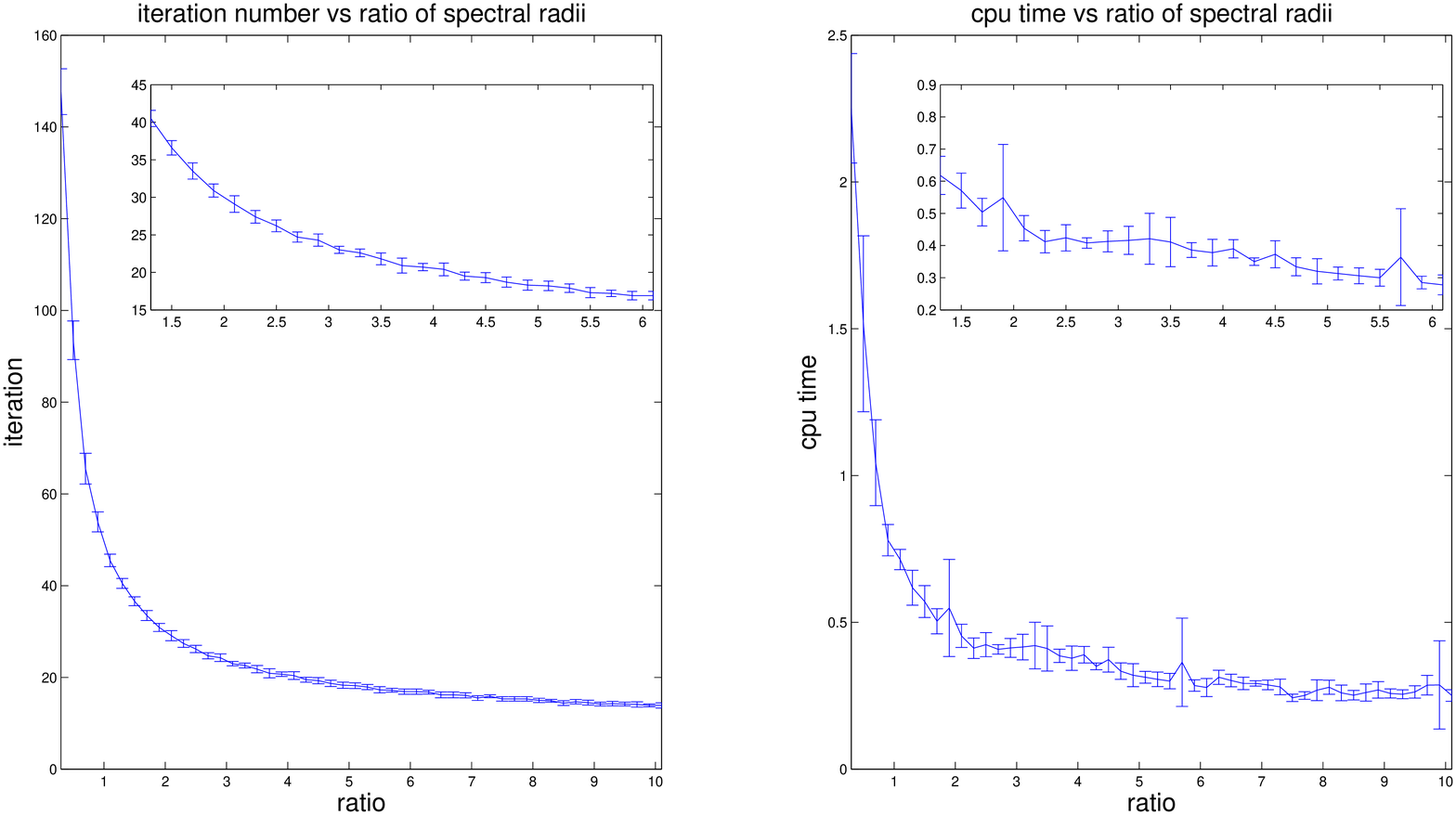}
\caption{The number of iterations (left), and
the CPU time (right) along with the ratios of the spectral radii (Case I).}
\end{figure}
\begin{figure}[htbp]
\centering
\includegraphics[width=6.8in]{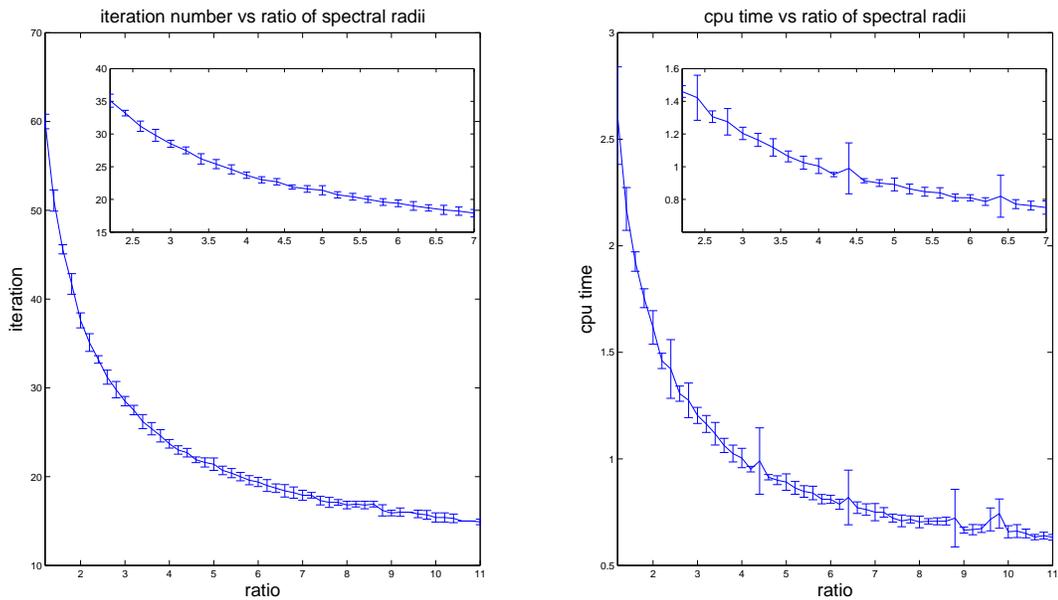}
\caption{The number of iterations (left), and
the CPU time (right) along with the ratios of the spectral radii (Case II).}
\end{figure}
\end{example}

 \begin{example}\label{exm:num3}
 This example tested tensors generated similar to Example~\ref{exm:num2}, with the only difference being that the partitions are
 \[
 |I_1|=\dots=|I_i|=2,\ \text{and }|I_{i+1}|=10\ \text{for }i\in [10],
 \]
 with again a unique genuine weakly irreducible block $I_{i+1}$.
The parameters were chosen as follows
\[
\text{Rt}= 2\ \text{and }\gamma=0.5\times 10^{-4}i\ \text{for all }i\in [10].
\]
 The information was recorded in Figure~4.  The magnitude is in logarithmic order.  We see that the computational effort increases exponentially along with the increasing of the number of blocks.
\begin{figure}[htbp]
\centering
\includegraphics[width=6.8in]{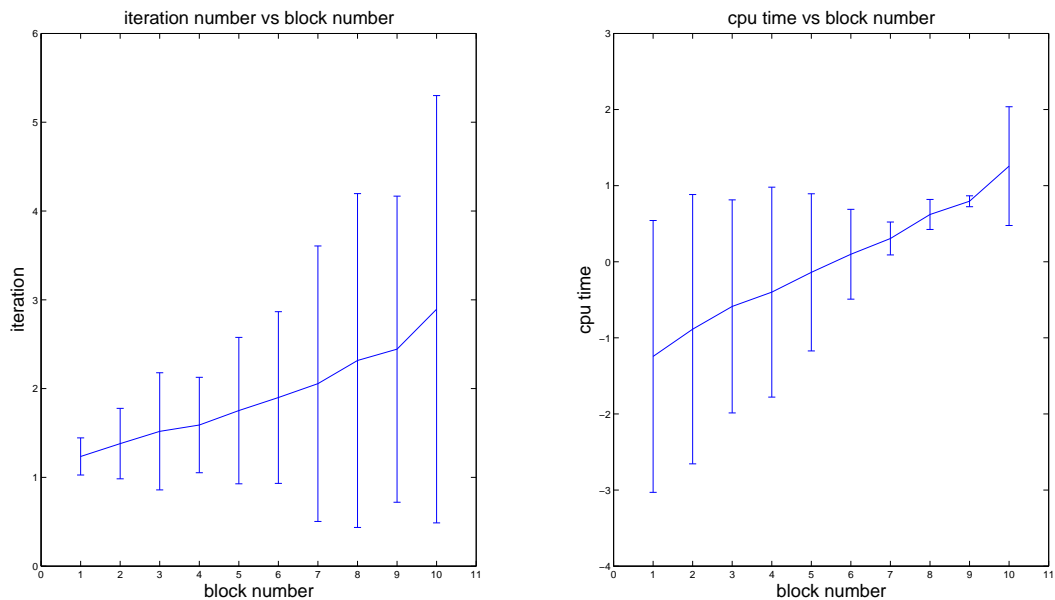}
\caption{The number of iterations (left), and
the CPU time (right) along with the numbers of blocks.}
\end{figure}
 \end{example}
 In conclusion, we see that Algorithm~\ref{algo:total} works very well. The choice for the initial point for the fixed iteration in Algorithm~\ref{algo:total} should affect the performance dramatically.

{\bf Acknowledgment}   We are thankful to Prof. Richard Brualdi and Prof. Jiayu Shao for their comments, and Mr. Weiyang Ding for reference \cite{DW}.

 }


\end{document}